\documentclass[12pt,a4paper]{article}

\usepackage[latin2]{inputenc}
\usepackage{amssymb}
\usepackage{paralist}
\usepackage{amsmath, calc}
\usepackage{amsfonts}
\usepackage{amsthm}
\usepackage{graphicx}
\usepackage{fullpage}
\usepackage{enumerate}
\newtheorem{theorem}{Theorem}
\newtheorem{lemma}{Lemma}
\newtheorem{conjecture}{Conjecture}

\newtheorem{problem}{Problem}

\begin{document}

\title{On the structure of sets which have coinciding representation functions}
\author{S\'andor Z. Kiss \thanks{Institute of Mathematics, Budapest
University of Technology and Economics, H-1529 B.O. Box, Hungary;
kisspest@cs.elte.hu;
This author was supported by the National Research, Development and Innovation Office NKFIH Grant No. K115288 and K109789, K129335. 
This paper was supported by the J\'anos Bolyai Research Scholarship of the Hungarian Academy of Sciences. Supported by the \'UNKP-18-4 New National Excellence Program of the Ministry of 
Human Capacities. Supported by the \'UNKP-19-4 New National Excellence Program 
of the Ministry for Innovation and Technology.}, Csaba
S\'andor \thanks{Institute of Mathematics, Budapest University of
Technology and Economics, H-1529 B.O. Box, Hungary, csandor@math.bme.hu.
This author was supported by the NKFIH Grants No. K109789, K129335. This paper was supported
by the J\'anos Bolyai Research Scholarship of the Hungarian Academy of Sciences.} 
}
\date{}
\maketitle

\begin{abstract}
\noindent For a set of nonnegative integers $A$, denote by $R_{A}(n)$ the number
of unordered representations of the integer $n$ as the sum of two
different terms from $A$. In this paper we partially describe the
structure of the sets, which have coinciding representation functions.


\end{abstract}

\section{Introduction}

Let $\mathbb{N}$ denote the set of nonnegative integers. For a given set $A \subseteq \mathbb{N}$, $A=\{a_1,a_2,\dots \}$ $(0\leq a_1<a_2< \dots)$, the additive representation functions $R_{h,A}^{(1)}(n)$,
$R_{h,A}^{(2)}(n)$ and $R_{h,A}^{(3)}(n)$ are defined in the following way:
\[
R_{h,A}^{(1)}(n)=|\{(a_{i_1},\dots ,a_{i_h}): a_{i_1}+\dots +a_{i_h}=n, a_{i_1},\dots ,a_{i_h} \in A\}|,
\]
\[
R_{h,A}^{(2)}(n)=|\{(a_{i_1},\dots ,a_{i_h}): a_{i_1}+\dots +a_{i_h}=n, a_{i_1}\leq a_{i_2}\leq \dots \leq a_{i_h}, a_{i_1},\dots, a_{i_h} \in A\}|,
\]
\[
R_{h,A}^{(3)}(n)=|\{(a_{i_1},\dots ,a_{i_h}): a_{i_1}+\dots +a_{i_h}=n, a_{i_1}<a_{i_2}<\dots <a_{i_h}, a_{i_1},\dots,a_{i_h} \in A\}|.
\]
For simplicity we write $R_{2,A}^{(3)}(n) = R_{A}(n)$. If $A$ is finite, let $|A|$ denote the cardinality of $A$.

The investigation of the partitions of the set of nonnegative integers with identical representation functions was a popular topic in the
 last few decades [1,3,4,5,7,9,11, 13,14]. It is easy to see that $R_{2,A}^{(1)}(n)$ is odd if and only if $\frac{n}{2}\in A$. It follows that for every positive integer $n$, $R_{2,C}^{(1)}(n)=R_{2,D}^{(1)}(n)$ holds if and only if $C=D$, where $C = \{c_{1}, c_{2}, \dots{}\}$ $(c_{1} < c_{2} < \dots{})$ and $D = \{d_{1}, d_{2}, \dots{}\}$ $(d_{1} < d_{2} < \dots{})$ are two sets of nonnegative integers. In [8], Nathanson gave a full description of the sets $C$ and $D$, which have identical representation functions $R_{2,C}^{(1)}(n)=R_{2,D}^{(1)}(n)$ from a certain point on. Namely, he proved the following theorem. Let $C(z) = \sum_{c \in C}z^{c}$, $D(z) = \sum_{d \in D}z^{d}$ be the generating functions of the sets $C$ and $D$, respectively.

\begin{theorem}
Let $C$ and $D$ be different infinite sets of nonnegative integers. Then $R_{2,C}^{(1)}(n)=R_{2,D}^{(1)}(n)$
holds from a certain point on if and only if there exist positive integers $n_{0}$,
$M$ and finite sets of nonnegative integers $F_{C}$, $F_{D}$, $T$ with
 $F_{C} \cup F_{D} \subset [0, Mn_{0} - 1]$, $T \subset [0, M - 1]$ such that
\[
C = F_{C} \cup \{lM + t: l \ge n_{0}, t \in T\},
\]
\[
D = F_{D} \cup \{lM + t: l \ge n_{0}, t \in T\},
\]
\[
1 - z^{M}|(F_{C}(z) - F_{D}(z))T(z).
\]
\end{theorem}

\noindent We conjecture in [6] that the above theorem of Nathanson can be generalized in the following way.

\begin{conjecture}
For $h > 2$ let $C$ and $D$ be different infinite sets of nonnegative
integers. Then $R_{h,C}^{(1)}(n)=R_{h,D}^{(1)}(n)$ holds
from a certain point on if and only if there exist positive integers $n_{0}$,
$M$ and finite sets $F_{C}$, $F_{D}$, $T$ with
 $F_{C} \cup F_{D} \subset [0, Mn_{0} - 1]$, $T \subset [0, M - 1]$ such that
\[
C = F_{C} \cup \{lM + t: l \ge n_{0}, t \in T\},
\]
\[
D = F_{D} \cup \{lM + t: l \ge n_{0}, t \in T\},
\]
\[
(1 - z^{M})^{h-1}|(F_{C}(z) - F_{D}(z))T(z)^{h-1}.
\]
\end{conjecture}

\noindent For $h = 3$, Kiss, Rozgonyi and S\'andor proved [6] Conjecture 1. In the general case when $h > 3$ we proved that if the conditions of Conjecture 1 hold then
$R_{h,C}^{(1)}(n)=R_{h,D}^{(1)}(n)$ holds from a certain point on.
Later, Rozgonyi and S\'andor [10] proved that the above conjecture holds when $h = p^{\alpha}$, where $\alpha \ge 1$ and $p$ is a prime.

It is easy to see that for any two different sets $C$, $D \subset \mathbb{N}$ we have  $R_{2,C}^{(2)}(n) \ne R_{2,D}^{(2)}(n)$ for some $n \in \mathbb{N}$.
Let $i$ denote the smallest index for which $c_i\neq d_i$, thus we may assume that $c_i<d_i$. It is
clear that $R_{2,C}^{(2)}(c_1+c_i)>R_{2,D}^{(2)}(c_1+c_i)$, which implies that there exists a nonnegative integer $n$ such that $R_{2,C}^{(2)}(n)\neq R_{2,D}^{(2)}(n)$.
We pose a problem about this representation function.

\begin{problem}
Determine all the sets of nonnegative integers $C$ and $D$ such that $R_{2,C}^{(2)}(n)=R_{2,D}^{(2)}(n)$ holds from a certain point on.
\end{problem}

\noindent In this paper, we focus on the representation function $R_{A}(n)$. We partially describe the structure of the sets, which have identical representation functions. To do this, 
we define the Hilbert cube which plays a crucial role in our results. Let $\{h_{1}, h_{2}, \dots{}\}$ $(h_1 < h_2 < \dots{})$ be finite or infinite set of positive integers. The set
\[
H(h_{1}, h_{2}, \dots{}) = \Big\{\sum_i\varepsilon_{i}h_{i}: \varepsilon_{i} \in \{0, 1\}\Big\}
\]
is called the Hilbert cube. The even part of a Hilbert cube is the set
\[
H_{0}(h_{1}, h_{2}, \dots{}) = \Big\{\sum_i\varepsilon_{i}h_{i}: \varepsilon_{i} \in \{0, 1\}, 2|\sum_i\varepsilon_{i}\Big\},
\]
and the odd part of a Hilbert cube is
\[
H_{1}(h_{1}, h_{2}, \dots{}) = \Big\{\sum_i\varepsilon_{i}h_{i}: \varepsilon_{i} \in \{0, 1\}, 2 \nmid \sum_i\varepsilon_{i}\Big\}.
\]
We say a Hilbert cube $H(h_{1}, h_{2}, \dots{})$ is half non-degenerated if the representation of any integer in $H_{0}(h_{1}, h_{2}, \dots{})$ and $H_{1}(h_{1}, h_{2}, \dots{})$ is unique, that is $\sum_i\varepsilon_{i}h_{i}\neq \sum_i\varepsilon_{i}^{'}h_{i}$ whenever $\sum_i\varepsilon_{i}\equiv \sum_i\varepsilon_{i}^{'} \hbox{ mod 2}$, where $\varepsilon_{i}^{'} \in \{0, 1\}$.

Many years ago, Selfridge and Straus [12] proved the following theorem about the cardinality of sets with identical representation functions. For the sake of completeness, we present the proof in Section 2.

\begin{theorem}
Let $C$ and $D$ be different finite sets of nonnegative integers such that for every positive integer $n$, $R_{C}(n) = R_{D}(n)$ holds. Then we have $|C| = |D| = 2^{l}$ for a nonnegative integer $l$.
\end{theorem}
\noindent If $0\in C$ and for $D=\{ d_1,d_2,\dots \}$, $0\leq d_1<d_2<\dots ,$ we have $R_C(m)=R_D(m)$ (sequences $C$ and $D$ are different), then $d_1>0$. Otherwise let us suppose that $c_i=d_i$ for $i=1,2,\dots ,n-1$, but $c_n<d_n$ which implies that $R_C(c_1+c_n)>R_D(c_1+c_n)$, a contradiction.

If $|C| = |D| = 1$ and $0 \in C$ with $R_{C}(n) = R_{D}(n)$, then we have $C=\{ 0\}$ and $D=\{ d_1\}$. Therefore, $C=H_0(d_1)$ and $D=H_1(d_1)$.

If $|C| = |D| = 2$ and $0 \in C$ with $R_{C}(n) = R_{D}(n)$, then $C=\{ 0,c_2\}$ and $D=\{ d_1,d_2\}$. In this case $1=R_C(0+c_2)$, and for $n \ne c_{2}$ we have $R_C(n) = 0$. Moreover, $1=R_C(d_1+d_2)$ and for $n \ne d_{1} + d_{2}$ we have $R_D(n) = 0$. This implies that $d_1+d_2=c_1+c_2=c_2$, that is $C=\{ 0,d_1+d_2\} =H_0(d_1,d_2)$ and $D=\{ d_1,d_2\}=H_1(d_1,d_2)$.

If $|C| = |D| = 4$ and $0 \in C$ with $R_{C}(n) = R_{D}(n)$, then let $C=\{ c_1,c_2,c_3,c_4\}$, $c_1=0$ and $D=\{ d_1,d_2,d_3,d_4\}$, where $d_{1} > 0$. Then we have
\[
c_1+c_2<c_1+c_3<c_1+c_4,c_2+c_3<c_2+c_4<c_3+c_4
\]
and
\[
d_1+d_2<d_1+d_3<d_1+d_4,d_2+d_3<d_2+d_4<d_3+d_4,
\]
which implies that $c_1+c_2=d_1+d_2$. Therefore, $c_2=d_1+d_2$ and $c_1+c_3=d_1+d_3$, thus we have $c_3=d_1+d_3$. If $c_2+c_3=d_2+d_3$, then $(d_1+d_2)+(d_1+d_3)=d_2+d_3$, that is $d_1=0$, a contradiction. Hence $c_2+c_3=d_1+d_4$, that is $(d_1+d_2)+(d_1+d_3)=d_1+d_4$. This implies that $d_4=d_1+d_2+d_3$. Finally $c_1+c_4=d_2+d_3$, that is $c_4=d_2+d_3$. Thus we have $C=\{ 0,d_1+d_2,d_1+d_3,d_2+d_3\}=H_0(d_1,d_2,d_3)$ and $D= \{d_{1}, d_{2}, d_{3}, d_{1}+d_{2}+d_{3}\} = H_1(d_1,d_2,d_3)$.
In the next step we prove that if the sets are even and odd parts of a Hilbert cube, then the corresponding representation functions are identical.
\begin{theorem}
Let $H(h_{1}, h_{2}, \dots{})$ be a half non-degenerated Hilbert cube.
If $C = H_{0}(h_{1}, h_{2}, \dots{})$ and $D = H_{1}(h_{1}, h_{2}, \dots{})$, then for every positive integer $n$,
$R_{C}(n) = R_{D}(n)$ holds.
\end{theorem}
\noindent It is easy to see that Theorem 3 is equivalent to Lemma 1 of Chen and Lev in [2]. First, they proved the finite case $H(h_1,\dots ,h_n)$ by induction on $n$, and the infinite case was a corollary of the finite case. For the sake of completeness, we give a different proof by using generating functions. Chen and Lev asked whether Theorem 3 described all different sets $C$ and $D$ of nonnegative integers such that $R_C(n)=R_D(n)$. The following conjecture is a simple generalization of the above question formulated by Chen and Lev but we use a different terminology.

\begin{conjecture}
Let $C$ and $D$ be different infinite sets of nonnegative integers with $0 \in C$. If for every positive integer $n$, $R_{C}(n) = R_{D}(n)$ holds, then there exist positive integers $d_{i_{1}}, d_{i_{2}}, \dots{} \in D$, where $d_{i_{1}} < d_{i_{2}} < \dots{},$ and a half non-degenerated Hilbert cube $H(d_{i_{1}}, d_{i_{2}}, \dots{})$ such that
\[
C = H_{0}(d_{i_{1}}, d_{i_{2}}, \dots{}),
\]
\[
D = H_{1}(d_{i_{1}}, d_{i_{2}}, \dots{}).
\]
\end{conjecture}
\noindent We showed above that Conjecture 2 is true for the finite case $l=0,1,2$. Unfortunately we could not settle the cases $l \ge 3$, which seems to be very complicated.
In Section 4 we prove the following weaker version of the above conjecture.

\begin{theorem}
Let $D = \{d_{1}, \dots{}, d_{2^{n}}\}$, $(0 < d_{1} < d_{2} < \dots{} < d_{2^{n}})$ be a set of nonnegative integers, where $d_{2^{k}+1} \ge 4d_{2^{k}}$, for $k = 0, \dots{}, n-1$ and $d_{2^{k}} \le d_{1} + d_{2} + d_{3} + d_{5} + \dots{} + d_{2^{i}+1} + \dots{} + d_{2^{k-1}+1}$ for $k = 2, \dots{}, n$. Let $C$ be a finite set of nonnegative integers such that $0 \in C$. If for every positive integer $m$, $R_{C}(m) = R_{D}(m)$ holds, then
\[
C = H_{0}(d_{1}, d_{2}, d_{3}, d_{5}, \dots{}, d_{2^{k}+1}, \dots ,d_{2^{n-1}+1}),
\]
and
\[
D = H_{1}(d_{1}, d_{2}, d_{3}, d_{5}, \dots{}, d_{2^{k}+1}, \dots ,d_{2^{n-1}+1}).
\]
\end{theorem}

\noindent For any sets of nonnegative integers $A$ and $B$ we define the sumset $A + B$ by
\[
A + B = \{a+b: a \in A, b \in B\}.
\]
The special case $b + A$ denotes the set $\{b + a: a \in A\}$, where $b$ is a fixed nonnegative integer.
Let $q\mathbb{N}$ denote the dilate of the set $\mathbb{N}$ by the factor $q$, that is, $q\mathbb{N}$ is the set of nonnegative integers divisible by $q$. Let $r_{A+B}(n)$ denote the number of solutions of the equation $a + b = n$, where $a \in A$, $b \in B$. In [2], Chen and Lev proved the following nice result.

\begin{theorem}
Let $l$ be a positive integer. Then there exist sets $C$ and $D$ of nonnegative integers such that $C\cup D=\mathbb{N}$, $C\cap D=2^{2l}-1+(2^{2l+1}-1)\mathbb{N}$ and for every positive integer $n$,
$R_C(n)=R_D(n)$ holds.
\end{theorem}

\noindent This theorem is an easy consequence of Theorem 3 by putting
\[
H(1,2,4,8,\dots ,2^{2l-1},2^{2l}-1,2^{2l+1}-1,2(2^{2l+1}-1),4(2^{2l+1}-1),8(2^{2l+1}-1),\dots ).
\]
The details can be found in the first part of the proof of Theorem 6. Chen and Lev [2] formulated the following conjecture.

\begin{conjecture}
Let $C$ and $D$ be different sets of nonnegative integers such that $C\cup D=\mathbb{N}$, $C\cap D=r+m\mathbb{N}$ with integers $r \ge 0$, $m \ge 2$. If for every positive integer
$n$, $R_C(n)=R_D(n)$ holds, then there exists an integer $l \ge 1$ such that $r=2^{2l}-1$ and $m=2^{2l+1}-1$.
\end{conjecture}

\noindent We formulate the following conjecture, which is a stronger version of the above conjecture of Chen and Lev.

\begin{conjecture}
Let $C$ and $D$ be different sets of nonnegative integers such that $C\cup D=\mathbb{N}$, $C\cap D=r+m\mathbb{N}$ with integers $r \ge 0$, $m \ge 2$. If for every positive integer $n$, $R_C(n)=R_D(n)$ holds, then there exists an integer $l \ge 1$ such that
\[
C=H_{0}(1,2,4,8,\dots, 2^{2l-1},2^{2l}-1,2^{2l+1}-1,2(2^{2l+1}-1),4(2^{2l+1}-1),8(2^{2l+1}-1),\dots ),
\]
and
\[
D=H_{1}(1,2,4,8,\dots, 2^{2l-1},2^{2l}-1,2^{2l+1}-1,2(2^{2l+1}-1),4(2^{2l+1}-1),8(2^{2l+1}-1),\dots ).
\]
\end{conjecture}

\noindent We prove that Conjecture 2 implies Conjecture 4.

\begin{theorem}
Assume that Conjecture 2 holds. Then there exist $C$ and $D$, different infinite sets of nonnegative integers, such that $C \cup D = \mathbb{N}$, $C\cap D = r + m\mathbb{N}$ with integers $r \ge 0$,
$m \ge 2$ and for every positive integer $n$, $R_C(n)=R_D(n)$ if and only if there exists an integer $l \ge 1$ such that
\[
C=H_{0}(1,2,4,8,\dots, 2^{2l-1},2^{2l}-1,2^{2l+1}-1,2(2^{2l+1}-1),4(2^{2l+1}-1),8(2^{2l+1}-1),\dots )
\]
and
\[
D=H_{1}(1,2,4,8,\dots, 2^{2l-1},2^{2l}-1,2^{2l+1}-1,2(2^{2l+1}-1),4(2^{2l+1}-1),8(2^{2l+1}-1),\dots ).
\]
\end{theorem}

\section{Proof of Theorem 2.}

\begin{proof}

By using the generating functions of the sets $C$ and $D$, we get that
\[
\sum_{n=1}^{\infty}R_{C}(n)z^{n} = \frac{C(z)^{2} - C(z^{2})}{2},
\]
\[
\sum_{n=1}^{\infty}R_{D}(n)z^{n} = \frac{D(z)^{2} - D(z^{2})}{2}.
\]
It follows that
\begin{equation}
R_C(n)=R_D(n) \hspace*{2mm} \textnormal{if and only if} \hspace*{2mm} C(z)^{2} - D(z)^{2} = C(z^{2}) - D(z^{2}).
\end{equation}
Let $l+1$ be the largest exponent of the factor $(z - 1)$ in $C(z) - D(z)$, i.e.,
\begin{equation}
C(z) - D(z) = (z - 1)^{l+1}p(z),
\end{equation}
where $p(z)$ is a polynomial and $p(1) \ne 0$. Substituting (2) back to (1) 
gives
\[
(C(z) + D(z))(z - 1)^{l+1}p(z) = (z^{2} - 1)^{l+1}p(z^{2}).
\]
Then we have $(C(z) + D(z))p(z) = (z + 1)^{l+1}p(z^{2})$. Substituting $z = 1$, we have $C(1) + D(1) = 2^{l+1}$, which implies that $|C| + |D| = 2^{l+1}$. On the other hand,
$$\binom{|C|}{2}=\sum_mR_C(m)=\sum_mR_D(m)=\binom{|D|}{2},$$ which gives $|C|=|D|$.
\end{proof}

\section{Proof of Theorem 3.}

\begin{proof}
By (1) we have to prove that $C(z)^{2} - D(z)^{2} = C(z^{2}) - D(z^{2})$. It is easy to see from the definition of $C$ and $D$ that
\[
\prod_i(1 - z^{h_{i}}) = \sum_{i_{1} < \dots{} < i_{t}}(-1)^{t}z^{h_{i_1} + \dots{} + h_{i_t}} = C(z) - D(z).
\]
On the other hand, clearly we have $C(z) + D(z) = \prod_i(1 + z^{h_{i}})$. Then we have
\[
C(z)^{2} - D(z)^{2} = \Big(C(z) - D(z)\Big)\Big(C(z) + D(z)\Big) = \prod_i(1 - z^{h_{i}}) \cdot \prod_i(1 + z^{h_{i}})
\]
\[
= \prod_i(1 - z^{2h_{i}}) = C(z^{2}) - D(z^{2}).
\]
The proof is completed.
\end{proof}

\section{Proof of Theorem 4.}

We apply induction on $n$. If $n = 0$, then $C = \{0\}$ and $D = \{d_{1}\}$. Therefore, for every positive integer $m$, we have $R_C(m)=R_D(m)=0$. If $n = 1$, then $C = \{0, c_{2}\}$ and $D = \{d_{1}, d_{2}\}$. Since $R_{C}(m) = R_{D}(m)$ for every positive integer $m$, it follows that $R_{D}(d_{1} + d_{2}) = 1 = R_{C}(d_{1} + d_{2})$. Then we have $C = \{0, d_{1} + d_{2}\} = H_{0}(d_{1}, d_{2})$ and $D = \{d_{1}, d_{2}\} = H_{1}(d_{1}, d_{2})$. Assume that the statement of Theorem 4 holds for $n = N - 1$. We will prove it for $n=N$. Let $D = \{d_{1}, \dots{}, d_{2^N}\}$ be a set of nonnegative integers, where $d_{2^{k}+1} \ge 4d_{2^{k}}$, for $k = 0, \dots{}, N - 1$ and $d_{2^{k}} \le d_{1} + d_{2} + d_{3} + d_{5} + \dots{} + d_{2^{i}+1} + \dots{} + d_{2^{k-1}+1}$ for $k = 2, \dots{}, N$. If $C$ is a set of nonnegative integers such that $0 \in C$ and for every positive integer $m$, $R_{C}(m)= R_{D}(m)$ holds, then we have to prove that
\[
C = H_{0}(d_{1}, d_{2}, d_{3}, d_{5}, \dots{},d_{2^k+1},\dots ,d_{2^{N-1}+1}),
\]
and
\[
D = H_{1}(d_{1}, d_{2}, d_{3}, d_{5}, \dots{},d_{2^k+1},\dots ,d_{2^{N-1}+1}).
\]
Define the sets
\[
C_{1} = \{c_{1}, \dots{}, c_{2^{N-1}}\}, \hspace*{20mm} C_{2} = C \setminus C_{1},
\]
and
\[
D_{1} = \{d_{1}, \dots{}, d_{2^{N-1}}\}, \hspace*{20mm} D_{2} = D \setminus D_{1}.
\]
We prove that for every positive integer $m$, we have
\begin{equation}
R_{C_{1}}(m) = R_{D_{1}}(m).
\end{equation}
Since $d_{2^{N-1}} \le \frac{1}{4}d_{2^{N-1}+1}$, it follows that for any $d_{i}$, $d_{j} \in D_1$ we have 
\[
d_{i} + d_{j} \le \frac{1}{4}d_{2^{N-1}+1} + \frac{1}{4}d_{2^{N-1}+1} = \frac{1}{2}d_{2^{N-1}+1}. 
\]
This implies that for every $\frac{1}{2}d_{2^{N-1}+1} \le m \le d_{2^{N-1}+1}$, we have $R_{D}(m) = 0$, which yields $R_{C}(m) = 0$. As $0 \in C$, we have a representation
$m = 0 + m$. It follows that $m \notin C$ for $\frac{1}{2}d_{2^{N-1}+1} \le m \le d_{2^{N-1}+1}$. We will show that 
\[
C_{1} = \Big[0, \frac{1}{3}d_{2^{N-1}+1}\Big[ \cap C, \hspace*{20mm} D_{1} = \Big[0, \frac{1}{3}d_{2^{N-1}+1}\Big[ \cap D.
\]
We distinguish two cases. In the first case we assume that $c_{2^{N-1}+1} \le \frac{d_{2^{N-1}+1}}{2}$. Then we have
\[
\binom{2^{N-1} + 1}{2} \le \sum_{m < d_{2^{N-1}+1}}R_{C}(m) = \sum_{m < d_{2^{N-1}+1}}R_{D}(m) = \binom{2^{N-1}}{2}
\]
which is a contradiction. In the second case we assume that $c_{2^{N-1}} > \frac{d_{2^{N-1}+1}}{2}$, which implies that $c_{2^{N-1}} \ge d_{2^{N-1}+1}$. Then we have
\[
\binom{2^{N-1} - 1}{2} \ge \sum_{m < d_{2^{N-1}+1}}R_{C}(m) = \sum_{m < d_{2^{N-1}+1}}R_{D}(m) = \binom{2^{N-1}}{2},
\]
which is impossible. Then we have $c_{2^{N-1}} \le \frac{1}{2}d_{2^{N-1}+1} < c_{2^{N-1}+1}$ and $d_{2^{N-1}+1} < c_{2^{N-1}+1}$, which implies that
\[
R_{C_{1}}(m) = \left\{
\begin{aligned}
0    \textnormal{, if } m \ge d_{2^{N-1}+1}\\
R_{C}(m) \textnormal{, if } m < d_{2^{N-1}+1}
\end{aligned} \hspace*{3mm},
\right.
\]
and
\[
R_{D_{1}}(m) = \left\{
\begin{aligned}
0    \textnormal{, if } m \ge d_{2^{N-1}+1}\\
R_{D}(m) \textnormal{, if } m < d_{2^{N-1}+1}
\end{aligned} \hspace*{3mm}.
\right.
\]
It follows that for every positive integer $m$, $R_{C_{1}}(m) = R_{D_{1}}(m)$, which proves (3). By the induction hypothesis we get that
\[
C_1=H_0(d_1,d_2,d_3,d_5,\dots, d_{2^k+1},\dots, ,d_{2^{N-2}+1})
\]
and
\[
D_1=H_1(d_1,d_2,d_3,d_5,\dots, d_{2^k+1},\dots, ,d_{2^{N-2}+1}).
\]
By Theorem 4, $d_{2^{k-1}+1} \le d_{2^{k}} \le \frac{1}{4}d_{2^{k}+1}$ for $1 \le k \le N-1$. This implies that $d_{2^{N-i}+1} \le \frac{1}{4^{i-1}}d_{2^{N-1}+1}$ for $i = 2, \dots{}, N$ and $d_{1} \le \frac{1}{4^{N}}d_{2^{N-1}+1}$. It follows that the largest element of the set $H(d_{1}, d_{2}, d_{3}, d_{5}, d_{9}, \dots{}, d_{2^{N-2}+1})$ is  
\begin{eqnarray*}
&&d_{1} + d_{2} + d_{3} + d_{5} + d_{9} + \dots{} + d_{2^{N-2}+1}\\ 
&\le& \frac{1}{4^{N}}d_{2^{N-1}+1} + \frac{1}{4^{N-1}}d_{2^{N-1}+1} + \dots{} + \frac{1}{4}d_{2^{N-1}+1} < \frac{1}{3}d_{2^{N-1}+1},
\end{eqnarray*}
 which implies that
\begin{equation}
C_{1} = \Big[0, \frac{1}{3}d_{2^{N-1}+1}\Big[ \cap C, \hspace*{20mm} D_{1} =
\Big]0, \frac{1}{3}d_{2^{N-1}+1}\Big[ \cap D.
\end{equation}
Therefore,
\begin{equation}
C_{1} + C_{1} \subset \Big[0, \frac{2}{3}d_{2^{N-1}+1}\Big[, \hspace*{20mm} D_{1} + D_{1} \subset
\Big[0, \frac{2}{3}d_{2^{N-1}+1}\Big[.
\end{equation}
Furthermore, for every $d \in D_{2}$, we have
\begin{eqnarray*}
&&d_{2^{N-1}+1} \le d \le d_{2^{N}} \le d_{2^{N-1}+1} + d_{2^{N-2}+1} + \dots{} + d_{2^{i}+1} + \dots{} + d_{2} + d_{1}\\
&\le& d_{2^{N-1}+1} + \frac{1}{4}d_{2^{N-1}+1} + \dots{} + \frac{1}{4^{N-i-1}}d_{2^{N-1}+1} + \dots{} \\ &+& \frac{1}{4^{N-1}}d_{2^{N-1}+1} + \frac{1}{4^{N}}d_{2^{N-1}+1}\\
&<& \frac{4}{3}d_{2^{N-1}+1}.
\end{eqnarray*}
Thus, we conclude
\begin{equation}
D_{1} + D_{2} \subset \Big[d_{2^{N-1}+1}, \frac{5}{3}d_{2^{N-1}+1}\Big[,
\end{equation}
and
\begin{equation}
D_{2} + D_{2} \subset \Big[2d_{2^{N-1}+1}, \frac{8}{3}d_{2^{N-1}+1}\Big[.
\end{equation}
It follows that
\begin{equation}
R_{C}(m) = 0 \hspace*{1mm} \text{for} \hspace*{1mm} m \ge \frac{8}{3}d_{2^{N-1}+1}.
\end{equation}
We now show that $c_{2^{N-1}+1} = d_{2^{N-1}+1} + d_{1}$. Assume that

\begin{equation}
c_{2^{N-1}+1} < d_{2^{N-1}+1} + d_{1}.
\end{equation}
Obviously, $c_{2^{N-1}+1} > d_{2^{N-1}+1}$. Since $c_{2^{N-1}+1} = c_{2^{N-1}+1} + 0$, then we have $1 \le R_{C}(c_{2^{N-1}+1}) = R_{D}(c_{2^{N-1}+1})$, which implies that
$c_{2^{N-1}+1} = d_{i} + d_{j}$, $i < j$, $d_{i}, d_{j} \in D$. If $j \le 2^{N-1}$, then by using the first condition in Theorem 4 we have
\[
c_{2^{N-1}+1} = d_{i} + d_{j} \le 2d_{2^{N-1}} \le \frac{1}{2}d_{2^{N-1}+1},
\]
which contradicts the inequality $c_{2^{N-1}+1} \geq d_{2^{N-1}+1}$. Moreover, when $j \ge 2^{N-1} + 1$, we have
\[
c_{2^{N-1}+1} = d_{i} + d_{j} \ge d_{1} + d_{2^{N-1}+1},
\]
which contradicts (9).

Assume that $c_{2^{N-1}+1} > d_{2^{N-1}+1} + d_{1}$.
Obviously, $1 \le R_{D}(d_{2^{N-1}+1} + d_{1}) = R_{C}(d_{2^{N-1}+1} + d_{1})$, which implies that $d_{1} + d_{2^{N-1}+1} = c_{i} + c_{j}$, $i < j$, $c_{i}, c_{j} \in C$. If $j \le 2^{N-1}$, then we have
\[
d_{1} + d_{2^{N-1}+1} = c_{i} + c_{j} \le 2c_{2^{N-1}} \le d_{2^{N-1}+1},
\]
which is impossible. Moreover, when $j \ge 2^{N-1} + 1$, we have
\[
d_{1} + d_{2^{N-1}+1} = c_{i} + c_{j} \ge c_{2^{N-1}+1} > d_{2^{N-1}+1} + d_{1}
\]
which is a contradiction. 

It follows that for every $c \in C$ with $c > c_{2^{N-1}+1}$, we have $c \le \frac{5}{3}d_{2^{N-1}+1}$. Otherwise,
$c + c_{2^{N-1}+1} \ge \frac{8}{3}d_{2^{N-1}+1}$ and then $R_{C}(c + c_{2^{N-1}+1}) \ge 1$, which contradicts (8).
By (4) and (8) we have
\begin{equation}
C_{1} + C_{2} \subset \Big[d_{2^{N-1}+1}, 2d_{2^{N-1}+1}\Big[,
\end{equation}
and
\begin{equation}
(C_{2} + C_{2}) \setminus \{2c_{2^{N}}\} \subset \Big[2d_{2^{N-1}+1}, \frac{8}{3}d_{2^{N-1}+1}\Big[.
\end{equation}

It sufficies to prove that 
\[
C_{2} = d_{2^{N-1}+1} + H_{1}(d_{1}, d_{2}, d_{3}, d_{5}, \dots{}, d_{2^{N-2}+1}) = d_{2^{N-1}+1} + D_{1}, 
\]
and 
\[
D_{2} = d_{2^{N-1}+1} + H_{0}(d_{1}, d_{2}, d_{3}, d_{5}, \dots{}, d_{2^{N-2}+1}) = d_{2^{N-1}+1} + C_{1}. 
\]
Define the sets
\[
C_{2,n} = \{c_{2^{N-1}+1}, c_{2^{N-1}+2}, \dots{}, c_{2^{N-1}+n}\},
\]
and
\[
D_{2,n} = \{d_{2^{N-1}+1}, d_{2^{N-1}+2}, \dots{}, d_{2^{N-1}+n}\}.
\]
Furthermore, define the sets
\[
C_{1} + C_{2,n} = \{p_{1}, p_{2}, \dots{}\}, \hspace*{1mm} (p_{1} < p_{2} < \dots{}), 
\]
\[
C_{2,n} + C_{2,n} = \{t_{1}, t_{2}, \dots{}\}, \hspace*{1mm} (t_{1} < t_{2} < \dots{}),
\]
and
\[
D_{1} + D_{2,n} = \{q_{1}, q_{2}, \dots{}\}, \hspace*{1mm} (q_{1} < q_{2} < \dots{}),
\]
\[
D_{2,n} + D_{2,n} = \{s_{1}, s_{2}, \dots{}\}, \hspace*{1mm} (s_{1} < s_{2} < \dots{}).
\]
Denote the first $2^{N-1}+n$ elements of the set
\[
H_{0}(d_1,d_2,d_3,d_5,\dots ,d_{{2^k}+1},\dots ,d_{2^{N-1}+1})
\]
by $H_{0}^{(n)}$, and let $H_{1}^{(n)}$ denote the first $2^{N-1}+n$ elements of the set
\[
H_{1}(d_1,d_2,d_3,d_5,\dots ,d_{{2^k}+1},\dots ,d_{2^{N-1}+1}).
\]
Now we prove by induction on $n$ that
\[
H^{(n)}_{0}=C_1\cup C_{2,n} \text{ and } H^{(n)}_{1} = D_1\cup D_{2,n}
\]
for $1 \le n \le 2^{N-1}$.
For $n = 1$ we have already proved that $D_{2,1}=\{ d_{2^{N-1}+1}\}$ and $C_{2,1}=\{ d_{2^{N-1}+1}+d_{1}\}$. It follows that $H^{(1)}_{0}=C_1\cup C_{2,1}$ and $H^{(1)}_{1} = D_1\cup D_{2,1}$.
Let us suppose that $H_{0}^{(n)}=C_1\cup C_{2,n}$ and  $H_{1}^{(n)}=D_1\cup D_{2,n}$ and we are going to prove that $H_{0}^{(n+1)}=C_1\cup C_{2,n+1}$ and $H_{1}^{(n+1)}=D_1\cup D_{2,n+1}$.
To prove $H_{0}^{(n+1)}=C_1\cup C_{2,n+1}$ and $H_{1}^{(n+1)}=D_1\cup D_{2,n+1}$, we need three lemmas. 
Let $i$ be the smallest index $u$ such that $r_{C_{1} + C_{2,n}}(p_{u}) > r_{D_{1} + D_{2,n}}(p_{u})$. If such an $i$ does not exist, then $p_{i} = +\infty$. Let $j$ be the smallest index $v$ such that $r_{C_{1} + C_{2,n}}(q_{v}) < r_{D_{1} + D_{2,n}}(q_{v})$. If such a $j$ does not exist, then $q_{j} = +\infty$. Let $k$ be the smallest index $w$ such that $R_{C_{2,n}}(t_{w})
> R_{D_{2,n}}(t_{w})$. If such a $k$ does not exist, then $t_{k} = +\infty$. Let $l$ be the smallest index $x$ such that $R_{C_{2,n}}(s_{x}) < R_{D_{2,n}}(s_{x})$. If such an $l$ does not exist, then $s_{l} = +\infty$. The following observations play a crucial role in the proof.

\begin{lemma}Let us suppose that $H_{0}^{(n)}=C_1\cup C_{2,n}$ and $H_{1}^{(n)}=D_1\cup D_{2,n}$. Then we have
\begin{itemize}
\item[(i)] $min\{p_{i}, c_{2^{N-1}+n+1}\} = min\{q_{j}, d_{1} + d_{2^{N-1}+n+1}\}$,
\item[(ii)] $min\{t_{k}, c_{2^{N-1}+1} + c_{2^{N-1}+n+1}\} = min\{s_{l}, d_{2^{N-1}+1} + d_{2^{N-1}+n+1}\}$.
\end{itemize}
\end{lemma}

\begin{proof} In the first step we prove (i).
We will prove that $p_i=+\infty$ is equivalent to $q_j=+\infty$ and for $p_i=q_j=+\infty$, we have $r_{C_{1} + C_{2,n}}(m)=r_{D_{1} + D_{2,n}}(m)$.
Assume that $p_{i} = +\infty$. Then by the definition of $p_i$, we have $r_{C_{1} + C_{2,n}}(p_{f}) \le r_{D_{1} + D_{2,n}}(p_{f})$ for every positive integer $f$. It follows that
\[
r_{C_{1} + C_{2,n}}(m) \le r_{D_{1} + D_{2,n}}(m)
\]
for every positive integer $m$.
On the other hand,
\[
2^{N-1}\cdot n = \sum_{m}r_{C_{1} + C_{2,n}}(m) \le \sum_{m}r_{D_{1} + D_{2,n}}(m) = 2^{N-1}\cdot n.
\]
Then $r_{C_{1} + C_{2,n}}(m) = r_{D_{1} + D_{2,n}}(m)$ for every positive integer $m$, which implies that
$q_{j} = +\infty$. Suppose that $q_{j} = +\infty$. Then by the definition of $q_j$, we have $r_{D_{1} + D_{2,n}}(q_{g}) \le
r_{C_{1} + C_{2,n}}(q_{g})$ for every positive integer $g$. It follows that
\[
r_{C_{1} + C_{2,n}}(m) \ge r_{D_{1} + D_{2,n}}(m)
\]
for every positive integer $m$.
Moreover, 
\[
2^{N-1}\cdot n = \sum_{m}r_{C_{1} + C_{2,n}}(m) \ge \sum_{m}r_{D_{1} + D_{2,n}}(m) = 2^{N-1}\cdot n.
\]
Then $r_{C_{1} + C_{2,n}}(m) = r_{D_{1} + D_{2,n}}(m)$ for every positive integer $m$, which implies that $p_{i} = +\infty$.

We distinguish two cases.

\noindent \textbf{Case 1.}
$p_{i} = +\infty$, $q_{j} = +\infty$, that is, $r_{C_{1} + C_{2,n}}(m)=r_{D_{1} + D_{2,n}}(m)$ for every positive integer $m$. Now we prove that $c_{2^{N-1}+n+1} = d_{1} + d_{2^{N-1}+n+1}$. Assume that  $c_{2^{N-1}+n+1} < d_{1} + d_{2^{N-1}+n+1}$.
Since $c_{2^{N-1}+n+1}=0+c_{2^{N-1}+n+1}$, where $0\in C_1$ but $c_{2^{N-1}+n+1}\in C_{2}\setminus C_{2,n}$ it follows from (5), (6), (7) and (10) that $R_{D}(c_{2^{N-1}+n+1}) = r_{D_{1} + D_{2,n}}(c_{2^{N-1}+n+1})$ and
$R_{C}(c_{2^{N-1}+n+1}) > r_{C_{1} + C_{2,n}}(c_{2^{N-1}+n+1})$. Then we have
\[
R_{D}(c_{2^{N-1}+n+1})=r_{D_1+D_{2,n}}(c_{2^{N-1}+n+1})=r_{C_1+C_{2,n}}(c_{2^{N-1}+n+1})<R_{C}(c_{2^{N-1}+n+1}),
\]
which is impossible.
Similarly, if $c_{2^{N-1}+n+1} > d_{1} + d_{2^{N-1}+n+1}$, then $R_{D}(d_{1} + d_{2^{N-1}+n+1}) > r_{D_{1} + D_{2,n}}(d_{1} + d_{2^{N-1}+n+1})$ because $d_1\in D_1$, $d_{2^{N-1}+n+1}\in D_{2}\setminus D_{2,n}$. It follows from (5), (6), (10) and (11) that $R_{C}(d_{1} + d_{2^{N-1}+n+1}) = r_{C_{1} + C_{2,n}}(d_{1} + d_{2^{N-1}+n+1})$. Then we have
\[
R_{C}(d_{1} + d_{2^{N-1}+n+1}) = r_{C_1+C_{2,n}}(d_{1} + d_{2^{N-1}+n+1}) = r_{D_1+D_{2,n}}(d_{1} + d_{2^{N-1}+n+1})
\]
\[
< R_{D}(d_{1} + d_{2^{N-1}+n+1}),
\]
which is a contradiction.

\noindent \textbf{Case 2.} $p_{i} < +\infty$ and $q_{j} < +\infty$. We have two subcases.

\textit{Case 2a.} $min\{p_{i}, c_{2^{N-1}+n+1}\} < min\{q_{j}, d_{1} + d_{2^{N-1}+n+1}\}$.

\noindent If $p_{i} \leq  c_{2^{N-1}+n+1}$, then obviously $p_{i} < d_{1} + d_{2^{N-1}+n+1}$, which implies by (5), (6), (7) and (10) that $R_{D}(p_{i}) = r_{D_{1} + D_{2,n}}(p_{i})$. By using the above facts and the definition of $p_{i}$, we get that
\[
R_{C}(p_{i}) \ge r_{C_{1} + C_{2,n}}(p_{i}) > r_{D_{1} + D_{2,n}}(p_{i}) = R_{D}(p_{i}),
\]
which contradicts the fact that $R_{C}(m) = R_{D}(m)$ for every positive integer $m$.
On the other hand, if $p_{i} > c_{2^{N-1}+n+1}$, then by the definition of $p_{i}$, $r_{C_{1} + C_{2,n}}(c_{2^{N-1}+n+1}) \le r_{D_{1} + D_{2,n}}(c_{2^{N-1}+n+1})$. Since $c_{2^{N-1}+n+1}=0+c_{2^{N-1}+n+1}$, $0\in C_1$ and $c_{2^{N-1}+n+1}\in C_{2}\setminus C_{2,n}$, we have
\begin{equation}
R_{C}(c_{2^{N-1}+n+1}) > r_{C_{1} + C_{2,n}}(c_{2^{N-1}+n+1}).
\end{equation}
The condition $min\{p_{i}, c_{2^{N-1}+n+1}\} < min\{q_{j}, d_{1} + d_{2^{N-1}+n+1}\}$ implies that $q_{j} > c_{2^{N-1}+n+1}$. It follows from the definition of $q_{j}$ that 
\[
r_{D_{1} + D_{2,n}}(c_{2^{N-1}+n+1}) \le r_{C_{1} + C_{2,n}}(c_{2^{N-1}+n+1}).
\]
We conclude that
\begin{equation}
r_{C_{1} + C_{2,n}}(c_{2^{N-1}+n+1}) = r_{D_{1} + D_{2,n}}(c_{2^{N-1}+n+1}).
\end{equation}
It follows from $0+c_{2^{N-1}+n+1} < d_{1} + d_{2^{N-1}+n+1}$, $0 \in C_{1}$, (5), (6), (7) and (10) that $r_{D_{1} + D_{2,n}}(c_{2^{N-1}+n+1}) = R_{D}(c_{2^{N-1}+n+1})$. Furthermore, we obtain from (12) and (13) that
\[
R_{D}(c_{2^{N-1}+n+1}) = r_{D_{1} + D_{2,n}}(c_{2^{N-1}+n+1}) = r_{C_{1} + C_{2,n}}(c_{2^{N-1}+n+1}) < R_{C}(c_{2^{N-1}+n+1}),
\]
which contradicts the fact that $R_{C}(m) = R_{D}(m)$ for every positive integer $m$.

\textit{Case 2b.} $min\{p_{i}, c_{2^{N-1}+n+1}\} > min\{q_{j}, d_{1} + d_{2^{N-1}+n+1}\}$.

\noindent If $q_{j} \leq d_1+d_{2^{N-1}+n+1} $, then obviously $q_{j} < c_{2^{N-1}+n+1}$, which implies from (5), (6), (10) and (11) that $R_{C}(q_{j}) = r_{C_{1} + C_{2,n}}(q_{j})$. By using the above facts and the definition of $q_{j}$, we get that
\[
R_{C}(q_{j}) = r_{C_{1} + C_{2,n}}(q_{j}) < r_{D_{1} + D_{2,n}}(q_{j}) \le R_{D}(q_{j}),
\]
which contradicts the fact that $R_{C}(m) = R_{D}(m)$ for every positive integer $m$.

\noindent Moreover, if $q_{j} > d_{1} + d_{2^{N-1}+n+1}$, then by the definition of $q_{j}$, 
\[
r_{C_{1} + C_{2,n}}(d_{1} + d_{2^{N-1}+n+1}) \ge r_{D_{1} + D_{2,n}}(d_{1} + d_{2^{N-1}+n+1}). 
\]
Since $d_{1} \in D_{1}$, 
$d_{2^{N-1}+n+1}\in D_{2}\setminus D_{2,n}$, we have
\begin{equation}
R_{D}(d_{1} + d_{2^{N-1}+n+1}) > r_{D_{1} + D_{2,n}}(d_{1} + d_{2^{N-1}+n+1}).
\end{equation}
The assumption $min\{p_{i}, c_{2^{N-1}+n+1}\} > min\{q_{j}, d_{1} + d_{2^{N-1}+n+1}\}$ implies that $p_{i} > d_{1} + d_{2^{N-1}+n+1}$. It follows from the definition of $p_{i}$ that 
\[
r_{C_{1} + C_{2,n}}(d_{1} + d_{2^{N-1}+n+1}) \le r_{D_{1} + D_{2,n}}(d_{1} + d_{2^{N-1}+n+1}). 
\]
We conclude that
\begin{equation}
r_{C_{1} + C_{2,n}}(d_{1} + d_{2^{N-1}+n+1}) = r_{D_{1} + D_{2,n}}(d_{1} + d_{2^{N-1}+n+1}).
\end{equation}
It follows from $min\{p_{i}, c_{2^{N-1}+n+1}\} > min\{q_{j}, d_{1} + d_{2^{N-1}+n+1}\}$ that $c_{2^{N-1}+n+1}>d_{1} + d_{2^{N-1}+n+1}$. Therefore, it follows from (5), (6), (10) and (11) that
\begin{equation}
r_{C_{1} + C_{2,n}}(d_{1} + d_{2^{N-1}+n+1}) = R_{C}(d_{1} + d_{2^{N-1}+n+1}).
\end{equation}
Furthermore, from (14),(15) and (16) we have
\[
R_{D}(d_{1} + d_{2^{N-1}+n+1}) > r_{D_{1} + D_{2,n}}(d_{1} + d_{2^{N-1}+n+1})
\]
\[
= r_{C_{1} + C_{2,n}}(d_{1} + d_{2^{N-1}+n+1}) = R_{C}(d_{1} + d_{2^{N-1}+n+1}),
\]
which contradicts the fact that $R_C(m)=R_D(m)$ for every positive integer $m$. The proof of (i) in Lemma 1 is completed.

The proof of (ii) in Lemma 1 is similar to the proof of (i). For the sake of completeness, we present it.
We prove that $s_l=+\infty$ is equivalent to $t_k=+\infty$ and in this case $R_{C_{2,n}}(m) =R_{D_{2,n}}(m)$ for every positive integer $m$.
If $t_{k} = +\infty$, then by the definition of $t_{k}$, we have $R_{C_{2,n}}(t_{f}) \le R_{D_{2,n}}(t_{f})$ for every positive integer $f$. Then
\[
R_{C_{2,n}}(m) \le R_{D_{2,n}}(m)
\]
for every positive integer $m$. On the other hand, we have
\[
\binom{n}{2} = \sum_{m}R_{C_{2,n}}(m) \le \sum_{m}R_{D_{2,n}}(m) = \binom{n}{2}.
\]
It follows that $R_{C_{2,n}}(m) = R_{D_{2,n}}(m)$ for every positive integer $m$, which implies that $s_{l} = +\infty$.
If $s_{l} = +\infty$, then by the definition of $s_{l}$, we have $R_{C_{2,n}}(s_{g}) \ge R_{D_{2,n}}(s_{g})$ for every positive integer $g$. Then
\[
R_{C_{2,n}}(m) \ge R_{D_{2,n}}(m)
\]
for every positive integer $m$. Furthermore,
\[
\binom{n}{2} = \sum_{m}R_{C_{2,n}}(m) \ge \sum_{m}R_{D_{2,n}}(m) = \binom{n}{2}.
\]
Then, we have $R_{C_{2,n}}(m) = R_{D_{2,n}}(m)$ for every positive integer $m$, which implies that $t_{k} = +\infty$.
We distinguish two cases.

\noindent \textbf{Case 1.} $t_{k} = +\infty$, $s_{l} = +\infty$, that is, $R_{C_{2,n}}(m) = R_{D_{2,n}}(m)$ for every positive integer $m$. Now, we prove that $d_{2^{N-1}+n+1} = d_{1} + c_{2^{N-1}+n+1}$. Assume that $d_{2^{N-1}+n+1} > d_{1} + c_{2^{N-1}+n+1}$. As $d_{1} + d_{2^{N-1}+1} + c_{2^{N-1}+n+1}=c_{2^{N-1}+1}+c_{2^{N-1}+n+1}$, where $c_{2^{N-1}+1}\in C_{2,n}$ and $c_{2^{N-1}+n+1}\in C_{2}\setminus C_{2,n}$, it follows that 
\[
R_{C}(c_{2^{N-1}+1} + c_{2^{N-1}+n+1}) > R_{C_{2,n}}(c_{2^{N-1}+1} + c_{2^{N-1}+n+1}). 
\]
On the other hand, we will show that 
\[
d_{2^{N-1}+1}+d_{2^{N-1}+n+1} > d_{1} + d_{2^{N-1}+1} +  c_{2^{N-1}+n+1} = c_{2^{N-1}+1} + c_{2^{N-1}+1} 
\]
and (5), (6), (7), (11) imply 
\[
R_{D}(c_{2^{N-1}+1} + c_{2^{N-1}+n+1}) = R_{D_{2,n}}(c_{2^{N-1}+1} + c_{2^{N-1}+n+1}). 
\]
It is clear from (11) that 
\[
c_{2^{N-1}+1} + c_{2^{N-1}+n+1} \in (C_{2} + C_{2}) \setminus \{2c_{2^{N}}\} \subset \Big[2d_{2^{N-1}+1}, \frac{8}{3}d_{2^{N-1}+1}\Big[. 
\]
By (5) and (6), we have $c_{2^{N-1}+1} + c_{2^{N-1}+n+1} \notin (D_{1} + D_{1}) \cup (D_{1} + D_{2})$. This implies that $R_{D}(c_{2^{N-1}+1} + c_{2^{N-1}+n+1}) = R_{D_{2}}(c_{2^{N-1}+1} + c_{2^{N-1}+n+1})$. Moreover, $D_{2} = D_{2,n} \cup (D_{2} \setminus D_{2,n})$, thus we have 
\[
R_{D_{2}}(m) = R_{D_{2,n}}(m) + 2r_{D_{2}+(D_{2} \setminus D_{2,n})}(m) + R_{D_{2} \setminus D_{2,n}}(m) 
\]
for any positive integer $m$. We conclude that for any positive integer $m$ with $2d_{2^{N-1}+1} \le m < d_{2^{N-1}+1} + d_{2^{N-1}+n+1}$, we have $r_{D_{2}+(D_{2} \setminus D_{2,n})}(m) = 0$ 
and 
\[
R_{D_{2} \setminus D_{2,n}}(m) = 0. 
\]
This implies that $R_{D_{2}}(m) = R_{D_{2,n}}(m)$, and then 
\[
R_{D}(c_{2^{N-1}+1} + c_{2^{N-1}+n+1}) = R_{D_{2,n}}(c_{2^{N-1}+1} + c_{2^{N-1}+n+1}). 
\]
Therefore,
\[
R_{D}(c_{2^{N-1}+1} + c_{2^{N-1}+n+1}) = R_{D_{2,n}}(c_{2^{N-1}+1} + c_{2^{N-1}+n+1})
\]
\[
= R_{C_{2,n}}(c_{2^{N-1}+1} + c_{2^{N-1}+n+1}) < R_{C}(c_{2^{N-1}+1} + c_{2^{N-1}+n+1}),
\]
which is impossible.
Similarly, if $d_{1} +c_{2^{N-1}+n+1} > d_{2^{N-1}+n+1}$, then it follows that $R_{D}(d_{2^{N-1}+1} + d_{2^{N-1}+n+1}) > R_{D_{2,n}}(d_{2^{N-1}+1} + d_{2^{N-1}+n+1})$ because $d_{2^{N-1}+1}\in D_{2,n}$ and $d_{2^{N-1}+n+1}\in D_{2}\setminus D_{2,n}$. It follows from (5), (7), (10), (11) and 
\[
d_{2^{N-1}+1} + d_{2^{N-1}+n+1} < d_{1} + d_{2^{N-1}+1} + c_{2^{N-1}+n+1} = c_{2^{N-1}+1} + c_{2^{N-1}+n+1} 
\]
that $R_{C_{2,n}}(d_{2^{N-1}+1} + d_{2^{N-1}+n+1}) = R_{C}(d_{2^{N-1}+1} + d_{2^{N-1}+n+1})$.
Then, we have
\[
R_C(d_{2^{N-1}+1} + d_{2^{N-1}+n+1}) = R_{C_{2,n}}(d_{2^{N-1}+1} + d_{2^{N-1}+n+1})
\]
\[
= R_{D_{2,n}}(d_{2^{N-1}+1} + d_{2^{N-1}+n+1}) < R_D(d_{2^{N-1}+1} + d_{2^{N-1}+n+1}),
\]
which is a contradiction.

\noindent \textbf{Case 2.} $t_{k} < +\infty$ and $s_{l} < +\infty$. We have two subcases.

\textit{Case 2a.} $min\{t_{k}, c_{2^{N-1}+1} + c_{2^{N-1}+n+1}\} < min\{s_{l}, d_{2^{N-1}+1} + d_{2^{N-1}+n+1}\}$. If $t_{k} \le c_{2^{N-1}+1} + c_{2^{N-1}+n+1}$, then obviously $t_{k} < d_{2^{N-1}+1} + d_{2^{N-1}+n+1}$. By using (5), (6), (7), (11) we have $R_{D}(t_{k}) = R_{D_{2,n}}(t_{k})$. By the definition of $t_{k}$, we get that
\[
R_{C}(t_{k}) \ge R_{C_{2,n}}(t_{k}) > R_{D_{2,n}}(t_{k}) = R_{D}(t_{k}),
\]
which contradicts the fact that $R_{C}(t_{k}) = R_{D}(t_{k})$.

\noindent On the other hand, if $t_{k} > c_{2^{N-1}+1} + c_{2^{N-1}+n+1}$, then by the definition of $t_{k}$, we have $R_{C_{2,n}}(c_{2^{N-1}+1} + c_{2^{N-1}+n+1})\leq R_{D_{2,n}}(c_{2^{N-1}+1} + c_{2^{N-1}+n+1})$. Moreover, it follows from $c_{2^{N-1}+1}\in C_{2,n}$ and $c_{2^{N-1}+n+1}\in C_{2}\setminus C_{2,n}$ that 
\[
R_{C}(c_{2^{N-1}+1} + c_{2^{N-1}+n+1}) > R_{C_{2,n}}(c_{2^{N-1}+1} + c_{2^{N-1}+n+1}). 
\]
Since $min\{t_{k}, c_{2^{N-1}+1} + c_{2^{N-1}+n+1}\} < min\{s_{l}, d_{2^{N-1}+1} + d_{2^{N-1}+n+1}\}$, we have $c_{2^{N-1}+1} + c_{2^{N-1}+n+1} <s_l$. By the definition of $s_{l}$, $R_{D_{2,n}}(c_{2^{N-1}+1} + c_{2^{N-1}+n+1}) \le R_{C_{2,n}}(c_{2^{N-1}+1} + c_{2^{N-1}+n+1})$. Then, we have 
\[
R_{C_{2,n}}(c_{2^{N-1}+1} + c_{2^{N-1}+n+1}) = R_{D_{2,n}}(c_{2^{N-1}+1} + c_{2^{N-1}+n+1}).
\]
By 
\[
c_{2^{N-1}+1} + c_{2^{N-1}+n+1} = min\{t_{k}, c_{2^{N-1}+1} + c_{2^{N-1}+n+1}\} 
\]
\[
< min\{s_{l}, d_{2^{N-1}+1} + d_{2^{N-1}+n+1}\} \le d_{2^{N-1}+1} + d_{2^{N-1} + n + 1} 
\]
and (5), (6), (7), (11), we have 
\[
R_{D_{2,n}}(c_{2^{N-1}+1} + c_{2^{N-1}+n+1}) = R_{D}(c_{2^{N-1}+1} + c_{2^{N-1}+n+1}), 
\]
and then
\[
R_{D}(c_{2^{N-1}+1} + c_{2^{N-1}+n+1}) = R_{D_{2,n}}(c_{2^{N-1}+1} + c_{2^{N-1}+n+1})
\]
\[
= R_{C_{2,n}}(c_{2^{N-1}+1} + c_{2^{N-1}+n+1}) < R_{C}(c_{2^{N-1}+1} + c_{2^{N-1}+n+1}),
\]
which contradicts the fact that 
\[
R_{C}(c_{2^{N-1}+1} + c_{2^{N-1}+n+1})= R_{D}(c_{2^{N-1}+1} + c_{2^{N-1}+n+1}).
\]

\textit{Case 2b.} $min\{t_{k}, c_{2^{N-1}+1} + c_{2^{N-1}+n+1}\} > min\{s_{l}, d_{2^{N-1}+1} + d_{2^{N-1}+n+1}\}$.

\noindent If $s_{l} \leq d_{2^{N-1}+1} + d_{2^{N-1}+n+1} $, then obviously $s_{l} < c_{2^{N-1}+1} + c_{2^{N-1}+n+1}$. By the definition of $s_{l}$ and (5), (7), (10), (11), it follows that
 $R_{C_{2,n}}(s_{l}) = R_{C}(s_{l})$, and then
\[
R_{C}(s_{l}) = R_{C_{2,n}}(s_{l}) < R_{D_{2,n}}(s_{l}) \le R_{D}(s_{l}),
\]
which contradicts the fact that $R_{C}(s_{l}) = R_{D}(s_{l})$.
On the other hand, if $s_{l} > d_{2^{N-1}+1} + d_{2^{N-1}+n+1}$, then by the definition of $s_{l}$, $R_{C_{2,n}}(d_{2^{N-1}+1} + d_{2^{N-1}+n+1}) \ge R_{D_{2,n}}(d_{2^{N-1}+1} + d_{2^{N-1}+n+1})$. Since 
$d_{2^{N-1}+1} \in D_{2}$, $d_{2^{N-1}+n+1} \in D_{2}\setminus D_{2,n}$, we have 
\[
R_{D}(d_{2^{N-1}+1} + d_{2^{N-1}+n+1}) > R_{D_{2,n}}(d_{2^{N-1}+1} + d_{2^{N-1}+n+1}).
\]
It follows from $min\{t_{k}, c_{2^{N-1}+1} + c_{2^{N-1}+n+1}\} > min\{s_{l}, d_{2^{N-1}+1} + d_{2^{N-1}+n+1}\}$ that $t_k>d_{2^{N-1}+1} + d_{2^{N-1}+n+1}$. 
By the definition of $t_{k}$, 
\[
R_{C_{2,n}}(d_{2^{N-1}+1} + d_{2^{N-1}+n+1}) \le R_{D_{2,n}}(d_{2^{N-1}+1} + d_{2^{N-1}+n+1}), 
\]
and then
$R_{C_{2,n}}(d_{2^{N-1}+1} + d_{2^{N-1}+n+1}) = R_{D_{2,n}}(d_{2^{N-1}+1} + d_{2^{N-1}+n+1})$.
By 
\[
d_{2^{N-1}+1} + d_{2^{N-1}+n+1} < c_{2^{N-1}+1} + c_{2^{N-1} + n + 1} 
\]
and (5), (7), (10), (11), we have
\[
R_{C_{2,n}}(d_{2^{N-1}+1} + d_{2^{N-1}+n+1}) = R_{C}(d_{2^{N-1}+1} + d_{2^{N-1}+n+1}), 
\]
and then
\[
R_{C}(d_{2^{N-1}+1} + d_{2^{N-1}+n+1}) = R_{C_{2,n}}(d_{2^{N-1}+1} + d_{2^{N-1}+n+1})
\]
\[
= R_{D_{2,n}}(d_{2^{N-1}+1} + d_{2^{N-1}+n+1})< R_{D}(d_{2^{N-1}+1} + d_{2^{N-1}+n+1}),
\]
which contradicts the fact that 
\[
R_C(d_{2^{N-1}+1} + d_{2^{N-1}+n+1})=R_D(d_{2^{N-1}+1} + d_{2^{N-1}+n+1}). 
\]
The proof of (ii) in Lemma 1 is completed.
\end{proof}

Let
\[
H = H(d_1,d_2,d_3,d_5,\dots, d_{2^k+1},\dots, ,d_{2^{N-1}+1})
\]
and
\[
H_0 = H_{0}(d_1,d_2,d_3,d_5,\dots, d_{2^k+1},\dots, ,d_{2^{N-1}+1})
\]
\[
H_1=H_{1}(d_1,d_2,d_3,d_5,\dots, d_{2^k+1},\dots, ,d_{2^{N-1}+1}).
\]
If $R_{H_{0}}(m) > 0$ or $R_{H_{1}}(m) > 0$, then
\[
m = \delta_{0}d_{1} + \sum_{i=1}^{N}\delta_{i}d_{2^{i-1}+1},
\]
where $\delta_{0}$, $\delta_{i} \in \{0, 1, 2\}$. It follows from $d_{2} \ge 4d_{1}$,
$d_{2^{k}+1} \ge 4d_{2^{k-1}+1}$, $(k = 1, \dots{}, N - 1)$ that when
\[
m^{'} = \delta_{0}^{'}d_{1} + \sum_{i=1}^{N}\delta_{i}^{'}d_{2^{i-1}+1},
\]
where $\delta_{0}^{'}$, $\delta_{i}^{'} \in \{0, 1, 2\}$ and $(\delta_{0}, \dots{}, \delta_{N}) \ne (\delta_{0}^{'}, \dots{}, \delta_{N}^{'})$, then $m \ne m^{'}$.
On the other hand, if 
\[
m = \delta_{0}d_{1} + \sum_{i=1}^{N}\delta_{i}d_{2^{i-1}+1},
\]
where $\delta_{0}$, $\delta_{i} \in \{0, 1, 2\}$, then $m = k + k^{'}$ with
\[
k = \varepsilon_{0}d_{1} + \sum_{i=1}^{N}\varepsilon_{i}d_{2^{i-1}+1},
\]
where $\varepsilon_{0}$, $\varepsilon_{i} \in \{0, 1\}$ and
\[
k^{'} = \varepsilon_{0}^{'}d_{1} + \sum_{i=1}^{N}\varepsilon_{i}^{'}d_{2^{i-1}+1},
\]
where $\varepsilon_{0}^{'}$, $\varepsilon_{i}^{'} \in \{0, 1\}$ if and only if
$\delta_{0} = \varepsilon_{0} + \varepsilon_{0}^{'}$ and
$\delta_{i} = \varepsilon_{i} + \varepsilon_{i}^{'}$, $1 \le i \le N$.

Let $H_{0,n}$ and $H_{1,n}$ denote the $2^{N-1}+n$th elements of $H_{0}$ and $H_{1}$,
respectively.
If
\[
H_{0,n} = \varepsilon_{0}d_{1} + \sum_{i=1}^{N}\varepsilon_{i}d_{2^{i-1}+1},
\]
where $\varepsilon_{0}$, $\varepsilon_{i} \in \{0, 1\}$, then it follows from $d_{2} \ge 4d_{1}$ and
$d_{2^{k}+1} \ge 4d_{2^{k-1}+1}$, $(k = 1, \dots{}, N - 1)$ that
\[
H_{1,n} = (1-\varepsilon_{0})d_{1} + \sum_{i=1}^{N}\varepsilon_{i}d_{2^{i-1}+1}.
\]
In the next step, we prove the following lemma.

\begin{lemma}
Let us suppose that $H_{0}^{(n)} = C_{1} \cup C_{2,n}$ and $H_{1}^{(n)} = D_{1} \cup D_{2,n}$ holds for some $1 \le n < 2^{N-1}$ . Let $H_{0,n+1} = \varepsilon_{0}d_{1} + \sum_{i=1}^{N}\varepsilon_{i}d_{2^{i-1}+1}$. If $\varepsilon_{0} = 0$  and $H_{0,n+1} = d_{2^{i_{1}-1}+1} + d_{2^{i_{2}-1}+1} + \dots{} + d_{2^{i_{t}-1}+1} + d_{2^{N-1}+1}$, where $1 \le i_{1} < i_{2} < \dots{} < i_{t} < N$, then we have
\begin{itemize}
\item[(i)] $q_{j} = H_{0,n+1}$,
\item[(ii)] $p_{i} > q_{j}$.
\end{itemize}
If $t > 1$, then 
\begin{itemize}
\item[(iii)] $s_{l} = 2d_{2^{i_{1}-1}+1} + d_{2^{i_{2}-1}+1} + \dots{} + d_{2^{i_{t}-1}+1} + 2d_{2^{N-1}+1}$, 
\item[(iv)] $t_{k} > s_{l}$.
\end{itemize}
If $t = 1$, then
\begin{itemize}
\item[(v)] $t_{k} = s_{l} = +\infty$.
\end{itemize}
\end{lemma}

\begin{proof}

\noindent We prove (i) and (ii) simultaneously. It is enough to show that
if $m < H_{0,n+1}$, then $r_{C_{1}+C_{2,n}}(m) = r_{D_{1}+D_{2,n}}(m)$ and
$r_{D_{1}+D_{2,n}}(H_{0,n+1}) > r_{C_{1}+C_{2,n}}(H_{0,n+1})$. If $m < d_{2^{N-1}+1}$, then it follows from (6) and (10) that $r_{C_{1}+C_{2,n}}(m) = r_{D_{1}+D_{2,n}}(m) = 0$.
If $d_{2^{N-1}+1} \le m < H_{0,n+1}$, then by using (5), (7), (11) and $H_{1,n+1} = H_{0,n+1} + d_{1}$, it follows that
$R_{H_{0}}(m) = r_{C_{1}+C_{2,n}}(m)$ and $R_{H_{1}}(m) = r_{D_{1}+D_{2,n}}(m)$. It follows from $R_{H_{0}}(m) = R_{H_{1}}(m)$ that $r_{C_{1}+C_{2,n}}(m) = r_{D_{1}+D_{2,n}}(m)$.
By using (5), (7), (11) and $H_{0,n+1} < H_{1,n+1}$, we get that $R_{H_{1}}(H_{0,n+1}) = r_{D_{1}+D_{2,n}}(H_{0,n+1})$. Since $H_{0,n+1} = 0 + H_{0,n+1}$, where
$0, H_{0,n+1} \in H_{0}$ and $H_{0,n+1} \notin C_{2,n}$, we have $R_{H_{0}}(H_{0,n+1}) > r_{C_{1}+C_{2,n}}(H_{0,n+1})$. It follows from $R_{H_{0}}(H_{0,n+1}) = R_{H_{1}}(H_{0,n+1})$ that $r_{D_{1}+D_{2,n}}(H_{0,n+1}) > R_{C_{1}+C_{2,n}}(H_{0,n+1})$, which proves (i) and (ii).

We prove (iii) and (iv) simultaneously. Let 
\[
M = 2d_{2^{i_{1}-1}+1} + d_{2^{i_{2}-1}+1} + \dots{} + d_{2^{i_{t}-1}+1} + 2d_{2^{N-1}+1}. 
\]
It is enough to show that
if $m < M$, then $R_{C_{2,n}}(m) = R_{D_{2,n}}(m)$ and $R_{C_{2,n}}(M) < R_{D_{2,n}}(M)$.
If $m < 2d_{2^{N-1}+1}$, then by using (7) and (11), we have $R_{C_{2,n}}(m) = R_{D_{2,n}}(m) = 0$. Let 
\[
2d_{2^{N-1}+1} \le m < d_{2^{i_{1}-1}+1} + d_{2^{i_{2}-1}+1} + \dots{} + d_{2^{i_{t}-1}+1} + 2d_{2^{N-1}+1}, 
\]
and write $m = h + h^{'}$ with $h, h^{'} \in H_{0}$.
By using (5) and (10), we get that $h, h^{'} \in H_{0} \setminus C_{1}$.
Since $h \ge d_{2^{N-1}+1}$, we have 
\[
h^{'} < d_{2^{i_{1}-1}+1} + d_{2^{i_{2}-1}+1} + \dots{} + d_{2^{i_{t}-1}+1} + d_{2^{N-1}+1} = H_{0,n+1}, 
\]
thus $h, h^{'} \in C_{2,n}$, which yields $R_{H_{0}}(m) = R_{C_{2,n}}(m)$. On the other hand, write $m = h + h^{'}$ with $h, h^{'} \in H_{1}$. By using (5) and (6), we get that $h, h^{'} \in H_{1} \setminus D_{1}$. Since $h \ge d_{2^{N-1}+1}$, we have 
\[
h^{'} < d_{2^{i_{1}-1}+1} + d_{2^{i_{2}-1}+1} + \dots{} + d_{2^{i_{t}-1}+1} + d_{2^{N-1}+1} = H_{0,n+1} < H_{1,n+1}, 
\]
thus $h, h^{'} \in D_{2,n}$, which yields $R_{H_{1}}(m) = R_{D_{2,n}}(m)$. It follows from
$R_{H_{0}}(m) = R_{H_{1}}(m)$ that $R_{C_{2,n}}(m) = R_{D_{2,n}}(m)$. Suppose that 
\[
d_{2^{i_{1}-1}+1} + d_{2^{i_{2}-1}+1} + \dots{} + d_{2^{i_{t}-1}+1} + 2d_{2^{N-1}+1}
\]
\[
\le m < 2d_{2^{i_{1}-1}+1} + d_{2^{i_{2}-1}+1} + \dots{} + d_{2^{i_{t}-1}+1} + 2d_{2^{N-1}+1}. 
\]
We can assume that
\[
m = \delta_{0}d_{1} + \sum_{i=1}^{u}d_{2^{x_{j}-1}+1} + \sum_{i=1}^{v}2d_{2^{y_{j}-1}+1}  
+ d_{2^{i_{1}-1}+1} + d_{2^{i_{2}-1}+1} + \dots{} + d_{2^{i_{t}-1}+1} + 2d_{2^{N-1}+1}, 
\]
where $\delta_{0} \in \{0, 1, 2\}$ and $1 \le x_{1} < x_{2} < \dots{} < x_{u} < i_{1}$ and $1 \le y_{1} < y_{2} < \dots{} < y_{v} < i_{1}$ and $x_{\alpha} \ne y_{\beta}$ are integers; otherwise, $R_{C_{2,n}}(m) = R_{D_{2,n}}(m) = 0$.   
Since $H_{0,n+1} = d_{2^{i_{1}-1}+1} + d_{2^{i_{2}-1}+1} + \dots{} + d_{2^{i_{t}-1}+1} + d_{2^{N-1}+1}$, then $t$ is odd, thus we can assume that $\delta_{0} + u + t$ is even;
otherwise, $R_{C_{2,n}}(m) = R_{D_{2,n}}(m) = 0$. We distinguish three cases.

\noindent \textbf{Case 1.} $\delta_{0} = 0$. Then, $u$ is odd. 
If $m = h + h^{'}$ with $h < h^{'}$ and $h, h^{'} \in H_{0}$, then it follows from (5) and (10) that $h, h^{'} \in H_{0} \setminus C_{1}$. It is clear that 
$h^{'} \notin C_{2,n}$ if and only if
\[
h^{'} = \sum_{j=1}^{w}d_{2^{z_{j}-1}+1} + \sum_{j=1}^{v}d_{2^{y_{j}-1}+1}  
+ d_{2^{i_{1}-1}+1} + d_{2^{i_{2}-1}+1} + \dots{} + d_{2^{i_{t}-1}+1} + d_{2^{N-1}+1},
\]
where $\{z_{1}, \dots{} ,z_{w}\} \subset \{x_{1}, \dots{} ,x_{u}\} \ne \emptyset$
and $w + v + t + 1$ is even. There are $2^{u-1}$ ways to choose the set $\{z_{1}, \dots{} ,z_{w}\}$, thus we have $R_{C_{2,n}}(m) = R_{H_{0}}(m) - 2^{u-1}$.

\noindent Furthermore, if $m = h + h^{'}$ with $h < h^{'}$ and $h, h^{'} \in H_{1}$, then it follows from (5) and (6) that $h, h^{'} \in H_{0} \setminus D_{1}$. It is clear that 
$h^{'} \notin D_{2,n}$ if and only if
\[
h^{'} = \sum_{j=1}^{w}d_{2^{z_{j}-1}+1} + \sum_{j=1}^{v}d_{2^{y_{j}-1}+1}  
+ d_{2^{i_{1}-1}+1} + d_{2^{i_{2}-1}+1} + \dots{} + d_{2^{i_{t}-1}+1} + d_{2^{N-1}+1},
\]
where $\{z_{1}, \dots{} ,z_{w}\} \subset \{x_{1}, \dots{} ,x_{u}\} \ne \emptyset$
and $w + v + t + 1$ is odd. There are $2^{u-1}$ ways to choose the set $\{z_{1}, \dots{} ,z_{w}\}$, thus we have $R_{C_{2,n}}(m) = R_{H_{1}}(m) - 2^{u-1}$.
As $R_{H_{0}}(m) = R_{H_{1}}(m)$, it follows that $R_{C_{2,n}}(m) = R_{D_{2,n}}(m)$.

\noindent \textbf{Case 2.} $\delta_{0} = 1$. Then, $u$ is even. 
If $m = h + h^{'}$ with $h < h^{'}$ and $h, h^{'} \in H_{0}$, then it follows from (5) and (10) that $h, h^{'} \in H_{0} \setminus C_{1}$. It is clear that 
$h^{'} \notin C_{2,n}$ if and only if
\[
h^{'} = \varepsilon_{0}d_{1} + \sum_{j=1}^{w}d_{2^{z_{j}-1}+1} + \sum_{j=1}^{v}d_{2^{y_{j}-1}+1}  
+ d_{2^{i_{1}-1}+1} + d_{2^{i_{2}-1}+1} + \dots{} + d_{2^{i_{t}-1}+1} + d_{2^{N-1}+1},
\]
where $\varepsilon_{0} \in \{0,1\}$ and $\{z_{1}, \dots{} ,z_{w}\} \subset \{x_{1}, \dots{} ,x_{u}\}$
and $\varepsilon_{0} + w + v + t + 1$ is even. If $u = 0$, then for a suitable
$\varepsilon_{0}$, there is only one possibility for $h^{'}$, thus we have $R_{C_{2,n}}(m) = R_{H_{0}}(m) - 1$. If $u > 0$, 
to choose the pairs $(\varepsilon_{0}, \{z_{1}, \dots{} ,z_{w}\})$, we 
have $2\cdot 2^{u-1} = 2^{u}$ possibilities, thus we have $R_{C_{2,n}}(m) = R_{H_{0}}(m) - 2^{u}$.
Moreover, if $m = h + h^{'}$, with $h < h^{'}$ and $h, h^{'} \in H_{1}$, then it follows from (5) and (6) that $h, h^{'} \in H_{1} \setminus D_{1}$. It is clear that 
$h^{'} \notin D_{2,n}$ if and only if
\[
h^{'} = \varepsilon_{0}d_{1} + \sum_{j=1}^{w}d_{2^{z_{j}-1}+1} + \sum_{j=1}^{v}d_{2^{y_{j}-1}+1}  
+ d_{2^{i_{1}-1}+1} + d_{2^{i_{2}-1}+1} + \dots{} + d_{2^{i_{t}-1}+1} + d_{2^{N-1}+1},
\]
where $\varepsilon_{0} \in \{0,1\}$ and $\{z_{1}, \dots{} ,z_{w}\} \subset \{x_{1}, \dots{} ,x_{u}\}$
and $\varepsilon_{0} + w + v + t + 1$ is odd. If $u = 0$, then for a suitable
$\varepsilon_{0}$, there is only one possibility for $h^{'}$ thus we have $R_{D_{2,n}}(m) = R_{H_{1}}(m) - 1$. When $u > 0$, 
to choose the pairs $(\varepsilon_{0}, \{z_{1}, \dots{} ,z_{w}\})$, we 
have $2\cdot 2^{u-1} = 2^{u}$ possibilities, thus we have $R_{D_{2,n}}(m) = R_{H_{1}}(m) - 2^{u}$. As $R_{H_{0}}(m) = R_{H_{1}}(m)$, it follows that $R_{C_{2,n}}(m) = R_{D_{2,n}}(m)$.

\noindent \textbf{Case 3.} $\delta_{0} = 2$. Then, $u$ is odd. 
If $m = h + h^{'}$ with $h < h^{'}$ and $h, h^{'} \in H_{0}$, then it follows from (5) and (10) that $h, h^{'} \in H_{0} \setminus C_{1}$. It is clear that 
$h^{'} \notin C_{2,n}$ if and only if
\[
h^{'} = d_{1} + \sum_{j=1}^{w}d_{2^{z_{j}-1}+1} + \sum_{j=1}^{v}d_{2^{y_{j}-1}+1}  
+ d_{2^{i_{1}-1}+1} + d_{2^{i_{2}-1}+1} + \dots{} + d_{2^{i_{t}-1}+1} + d_{2^{N-1}+1},
\]
where $\{z_{1}, \dots{} ,z_{w}\} \subset \{x_{1}, \dots{} ,x_{u}\} \ne \emptyset$
and $1 + w + v + t + 1$ is even. There are $2^{u-1}$ ways to choose the set $\{z_{1}, \dots{} ,z_{w}\}$, thus we have $R_{C_{2,n}}(m) = R_{H_{0}}(m) - 2^{u-1}$.

\noindent On the other hand, if $m = h + h^{'}$ with $h < h^{'}$ and $h, h^{'} \in H_{1}$, then it follows from (5) and (6) that $h, h^{'} \in H_{0} \setminus D_{1}$. It is clear that 
$h^{'} \notin D_{2,n}$ if and only if
\[
h^{'} = d_{1} + \sum_{j=1}^{w}d_{2^{z_{j}-1}+1} + \sum_{j=1}^{v}d_{2^{y_{j}-1}+1}  
+ d_{2^{i_{1}-1}+1} + d_{2^{i_{2}-1}+1} + \dots{} + d_{2^{i_{t}-1}+1} + d_{2^{N-1}+1},
\]
where $\{z_{1}, \dots{} ,z_{w}\} \subset \{x_{1}, \dots{} ,x_{u}\} \ne \emptyset$
and $1 + w + v + t + 1$ is odd. There are $2^{u-1}$ ways to choose the set $\{z_{1}, \dots{} ,z_{w}\}$, thus we have $R_{C_{2,n}}(m) = R_{H_{1}}(m) - 2^{u-1}$. As $R_{H_{0}}(m) = R_{H_{1}}(m)$, 
it follows that $R_{C_{2,n}}(m) = R_{D_{2,n}}(m)$. Let 
\[
M = 2d_{2^{i_{1}-1}+1} + d_{2^{i_{2}-1}+1} + \dots{} + d_{2^{i_{t}-1}+1} + 2d_{2^{N-1}+1}.
\]
Now, we prove $R_{D_{2,n}}(M) = R_{H_{1}}(M)$. Assume that 
\[
M = 2d_{2^{i_{1}-1}+1} + d_{2^{i_{2}-1}+1} + \dots{} + d_{2^{i_{t}-1}+1} + 2d_{2^{N-1}+1} = h + h^{'},
\]
where $h, h^{'} \in H_{1}$ with $h < h^{'}$. Then, it follows from (5) and (6) that $h, h^{'} \in H_{1} \setminus D_{1}$. It follows that
\[
h^{'} = d_{2^{i_{1}-1}+1} + d_{2^{z_{1}-1}+1} + \dots{} + d_{2^{z_{w}-1}+1} + d_{2^{N-1}},
\]
where $\{z_{1}, \dots{} ,z_{w}\} \subset \{i_{2}, \dots{} ,i_{t}\}$. Thus, we have 
\[
h^{'} \le d_{2^{i_{1}-1}+1} + d_{2^{i_{2}-1}+1} + \dots{} + d_{2^{i_{t}-1}+1} + d_{2^{N-1}+1} = H_{0,n+1} < H_{1,n+1}, 
\]
which implies 
that $R_{D_{2,n}}(M) = R_{H_{1}}(M)$. On the other hand,
\[
M = (d_{2^{i_{1}-1}+1} + d_{2^{N-1}+1}) + (d_{2^{i_{1}-1}+1} + d_{2^{i_{2}-1}+1} + \dots{} + d_{2^{i_{t}-1}+1} + d_{2^{N-1}+1}),
\]
where $d_{2^{i_{1}-1}+1} + d_{2^{N-1}+1} \in C_{2,n}$ and 
\[
d_{2^{i_{1}-1}+1} + d_{2^{i_{2}-1}+1} + \dots{} + d_{2^{i_{t}-1}+1} + d_{2^{N-1}+1} \notin C_{2,n}, 
\]
and we have 
\[
d_{2^{i_{1}-1}+1} + d_{2^{N-1}+1} < d_{2^{i_{1}-1}+1} + d_{2^{i_{2}-1}+1} + \dots{} + d_{2^{i_{t}-1}+1} + d_{2^{N-1}+1} 
\]
because $t \ge 2$. This gives $R_{H_{0}}(M) > R_{C_{2,n}}(M)$. It follows from $R_{H_{0}}(M) = R_{H_{1}}(M)$
that $R_{D_{2,n}}(M) > R_{C_{2,n}}(M)$. 
Assume that $t = 1$, that is, $H_{0,n+1} = d_{2^{i_{1}-1}+1} + d_{2^{N-1}+1}$. The previous argument shows that for $m < 2d_{2^{i_{1}-1}+1} + 2d_{2^{N-1}+1}$,
we have $R_{C_{2,n}}(m) = R_{D_{2,n}}(m)$. Moreover, if $m \ge 2d_{2^{i_{1}-1}} + 2d_{2^{N-1}+1} = 2H_{0,n+1}$ and $R_{C_{2,n}}(m) \ne 0$ or $R_{D_{2,n}}(m) \ne 0$, then 
\[
m = \delta_{0}d_{1} + \sum_{u=1}^{s}\delta_{u}d_{2^{j_{u}-1}+1} + 2d_{2^{N-1}+1},
\]
where $\delta_{0} \in \{0,1,2\}$, $\delta_{u} \in \{1,2\}$, $1 \le j_{1} < j_{2} < \dots{} < j_{s} < N$ and $j_{s} \ge i_{1}$. If $m = h + h^{'}$ with
$h, h^{'} \in H_{0}$ or $h, h^{'} \in H_{1}$, $h < h^{'}$, then
\[
h^{'} = \varepsilon_{0}d_{1} + \sum_{l=1}^{r}d_{2^{h_{l}-1}+1} + d_{2^{N-1}+1},
\]
where $1 \le h_{1} < j_{2} < \dots{} < h_{r}$ and $h_{r} = j_{s} \ge i_{1}$. Hence, we have $h^{'} \ge  H_{0,n+1}$. It follows that 
$h^{'} \notin C_{2,n}$ and $h^{'} \notin D_{2,n}$, which implies that $R_{C_{2,n}}(m) = R_{D_{2,n}}(m) = 0$. This proves $s_{l} = t_{k} = +\infty$.
\end{proof}

\begin{lemma}
For $1 \le n < 2^{N-1}$, let $H_{0}^{(n)} = C_{1} \cup C_{2,n}$ and $H_{1}^{(n)} = D_{1} \cup D_{2,n}$. Let $H_{0,n+1} = \varepsilon_{0}d_{1} + \sum_{i=1}^{N}\varepsilon_{i}d_{2^{i-1}+1}$. If $\varepsilon_{0} = 1$ and $H_{0,n+1} = d_{1} + d_{2^{i_{1}-1}+1} + d_{2^{i_{2}-1}+1} + \dots{} + d_{2^{i_{t}-1}+1} + d_{2^{N-1}+1}$, where $1 \le i_{1} < i_{2} < \dots{} < i_{t} < N$, then we have
\begin{itemize}
\item[(i)] $p_{i} = 2d_{2^{i_{1}-1}+1} + d_{2^{i_{2}-1}+1} + \dots{} + d_{2^{i_{t}-1}+1} + d_{2^{N-1}+1}$,
\item[(ii)] $q_{j} > p_{i}$,
\item[(iii)] $t_{k} = d_{2^{i_{1}-1}+1} + d_{2^{i_{2}-1}+1} + \dots{} + d_{2^{i_{t}-1}+1} + 2d_{2^{N-1}+1}$, 
\item[(iv)] $s_{l} > t_{k}$.
\end{itemize}
\end{lemma}

\begin{proof}
\noindent We prove (i) and (ii) simultaneously. It is enough to show that
for $\varepsilon_{0} = 1$ and 
\[
H_{0,n+1} = d_{1} + d_{2^{i_{1}-1}+1} + d_{2^{i_{2}-1}+1} + \dots{} + d_{2^{i_{t}-1}+1} + d_{2^{N-1}+1},
\]
if $m < 2d_{2^{i_{1}-1}+1} + d_{2^{i_{2}-1}+1} + \dots{} + d_{2^{i_{t}-1}+1} + d_{2^{N-1}+1}$, 
then $r_{C_{1}+C_{2,n}}(m) = r_{D_{1}+D_{2,n}}(m)$ and
$r_{D_{1}+D_{2,n}}(K) < r_{C_{1}+C_{2,n}}(K)$, where 
\[
K =  2d_{2^{i_{1}-1}+1} + d_{2^{i_{2}-1}+1} + \dots{} + d_{2^{i_{t}-1}+1} + d_{2^{N-1}+1}. 
\]
If $m < d_{2^{N-1}+1}$, then it follows from (6) and (10) that $r_{C_{1}+C_{2,n}}(m) = r_{D_{1}+D_{2,n}}(m) = 0$.
Assume that 
\[
d_{2^{N-1}+1} \le m <  d_{2^{i_{1}-1}+1} + d_{2^{i_{2}-1}+1} + \dots{} + d_{2^{i_{t}-1}+1} + d_{2^{N-1}+1}. 
\]
If $m = h + h^{'}$ with $h < h^{'}$ and $h, h^{'} \in H_{0}$, then
it follows from (5) and (11) that $h \in C_{1}$ and $h^{'} \in H_{0} \setminus C_{1}$.
Since 
\[
h^{'} < d_{2^{i_{1}-1}+1} + d_{2^{i_{2}-1}+1} + \dots{} + d_{2^{i_{t}-1}+1} + d_{2^{N-1}+1} < H_{0,n+1}, 
\]
we have $h^{'} \in C_{2,n}$, which yields $R_{H_{0}}(m) = r_{C_{1}+C_{2,n}}(m)$. 
Upon writing $m = h + h^{'}$ with $h < h^{'}$ and $h, h^{'} \in H_{1}$, it follows from (5) and (7) that $h \in D_{1}$ and $h^{'} \in H_{1} \setminus D_{1}$.
Since 
\[
h^{'} < d_{2^{i_{1}-1}+1} + d_{2^{i_{2}-1}+1} + \dots{} + d_{2^{i_{t}-1}+1} + d_{2^{N-1}+1} = H_{1,n+1}, 
\]
we have $h^{'} \in D_{2,n}$. As $R_{H_{0}}(m) = R_{H_{0}}(m)$, we have $r_{D_{1}+D_{2,n}}(m) = r_{C_{1}+C_{2,n}}(m)$.
Suppose that 
\[
d_{2^{i_{1}-1}+1} + d_{2^{i_{2}-1}+1} + \dots{} + d_{2^{i_{t}-1}+1} + d_{2^{N-1}+1} 
\]
\[
\le m < 2d_{2^{i_{1}-1}+1} + d_{2^{i_{2}-1}+1} + \dots{} + d_{2^{i_{t}-1}+1} + d_{2^{N-1}+1}. 
\]
Then, we may assume that $m$ can be written in the form
\[
m = \delta_{0}d_{1} + \sum_{i=1}^{u}d_{2^{x_{j}-1}+1} + \sum_{i=1}^{v}2d_{2^{y_{j}-1}+1}  
+ d_{2^{i_{1}-1}+1} + d_{2^{i_{2}-1}+1} + \dots{} + d_{2^{i_{t}-1}+1} + d_{2^{N-1}+1}, 
\]
where $\delta_{0} \in \{0, 1, 2\}$ and $1 \le x_{1} < x_{2} < \dots{} < x_{u} < i_{1}$ and $1 \le y_{1} < y_{2} < \dots{} < y_{v} < i_{1}$, and $x_{\alpha} \ne y_{\beta}$ are integers; otherwise, 
$r_{C_{1}+C_{2,n}}(m) = r_{D_{1}+D_{2,n}}(m) = 0$.   

Since $H_{0,n+1} = d_{1} + d_{2^{i_{1}-1}+1} + d_{2^{i_{2}-1}+1} + \dots{} + d_{2^{i_{t}-1}+1} + d_{2^{N-1}+1}$, we have $t$ is even, thus 
$\delta_{0} + u + 2v + t + 1$ is even, which implies that $\delta_{0} + u$ is odd.
We distinguish three cases.

\noindent \textbf{Case 1.} $\delta_{0} = 0$. Then, $u$ is odd.  
If $m = h + h^{'}$ with $h < h^{'}$ and $h, h^{'} \in H_{0}$, then it follows from (5) and (10) that $h \in C_{1}$ and $h^{'} \in H_{0} \setminus C_{1}$. 
It is clear that 
$h^{'} \notin C_{2,n}$ if and only if $h^{'}$ can be written in the form 
\[
h^{'} = \sum_{j=1}^{w}d_{2^{z_{j}-1}+1} + \sum_{j=1}^{v}d_{2^{y_{j}-1}+1}  
+ d_{2^{i_{1}-1}+1} + d_{2^{i_{2}-1}+1} + \dots{} + d_{2^{i_{t}-1}+1} + d_{2^{N-1}+1},
\]
where $\{z_{1}, \dots{} ,z_{w}\} \subset \{x_{1}, \dots{} ,x_{u}\} \ne \emptyset$
and $w + v + t + 1$ is even. There are $2^{u-1}$ ways to choose the set $\{z_{1}, \dots{} ,z_{w}\}$, thus we have $r_{C_{1}+C_{2,n}}(m) = R_{H_{0}}(m) - 2^{u-1}$.
Furthermore, if $m = h + h^{'}$ with $h < h^{'}$ and $h, h^{'} \in H_{1}$, then it follows from (5) and (7) that $h \in D_{1}$ and $h^{'} \in H_{1} \setminus D_{1}$. It is clear that 
$h^{'} \notin D_{2,n}$ if and only if $h^{'}$ can be written in the form
\[
h^{'} = \sum_{j=1}^{w}d_{2^{z_{j}-1}+1} + \sum_{j=1}^{v}d_{2^{y_{j}-1}+1}  
+ d_{2^{i_{1}-1}+1} + d_{2^{i_{2}-1}+1} + \dots{} + d_{2^{i_{t}-1}+1} + d_{2^{N-1}+1},
\]
where $\{z_{1}, \dots{} ,z_{w}\} \subset \{x_{1}, \dots{} ,x_{u}\} \ne \emptyset$
and $w + v + t + 1$ is odd. There are $2^{u-1}$ ways to choose the set $\{z_{1}, \dots{} ,z_{w}\}$, thus we have $r_{D_{1}+D_{2,n}}(m) = R_{H_{1}}(m) - 2^{u-1}$.
As $R_{H_{0}}(m) = R_{H_{1}}(m)$, it follows that $r_{C_{1}+C_{2,n}}(m) = r_{D_{1}+D_{2,n}}(m)$.

\noindent \textbf{Case 2.} $\delta_{0} = 1$. Then, $1 + u + 2v + t + 1$
is even, which implies that $u$ is even.
  
If $m = h + h^{'}$ with $h < h^{'}$ and $h, h^{'} \in H_{0}$, then $h \in C_{1}$ and $h^{'} \in H_{0} \setminus C_{1}$. It is clear that 
$h^{'} \notin C_{2,n}$ if and only if $h^{'}$ can be written in the form
\[
h^{'} = \varepsilon_{0}d_{1} + \sum_{j=1}^{w}d_{2^{z_{j}-1}+1} + \sum_{j=1}^{v}d_{2^{y_{j}-1}+1}  
+ d_{2^{i_{1}-1}+1} + d_{2^{i_{2}-1}+1} + \dots{} + d_{2^{i_{t}-1}+1} + d_{2^{N-1}+1},
\]
where $\varepsilon_{0} \in \{0,1\}$ and $\{z_{1}, \dots{} ,z_{w}\} \subset \{x_{1}, \dots{} ,x_{u}\}$
and $\varepsilon_{0} + w + v + t + 1$ is even. If $u = 0$, then $\{z_{1}, \dots{} ,z_{w}\} = \emptyset$ and for a suitable
$\varepsilon_{0}$, there is only one way to choose $h^{'}$, and so $r_{C_{1}+C_{2,n}}(m) = R_{H_{0}}(m) - 1$. When $u > 0$ is even, 
to choose the pairs $(\varepsilon_{0}, \{z_{1}, \dots{} ,z_{w}\})$ we 
have $2\cdot 2^{u-1} = 2^{u}$ possibilities, thus we have $r_{C_{2,n}}(m) = R_{H_{0}}(m) - 2^{u}$.

On the other hand, if $m = h + h^{'}$ with $h < h^{'}$ and $h, h^{'} \in H_{0}$, then $h \in D_{1}$ and $h^{'} \in H_{1} \setminus D_{1}$. It is clear that 
$h^{'} \notin D_{2,n}$ if and only if $h^{'}$ can be written in the form
\[
h^{'} = \varepsilon_{0}d_{1} + \sum_{j=1}^{w}d_{2^{z_{j}-1}+1} + \sum_{j=1}^{v}d_{2^{y_{j}-1}+1}  
+ d_{2^{i_{1}-1}+1} + d_{2^{i_{2}-1}+1} + \dots{} + d_{2^{i_{t}-1}+1} + d_{2^{N-1}+1},
\]
where $\varepsilon_{0} \in \{0,1\}$ and $\{z_{1}, \dots{} ,z_{w}\} \subset \{x_{1}, \dots{} ,x_{u}\}$
and $\varepsilon_{0} + w + v + t + 1$ is odd. If $u = 0$, then $\{z_{1}, \dots{} ,z_{w}\} = \emptyset$, and for a suitable
$\varepsilon_{0}$, there is only one way to choose $h^{'}$, and so $r_{D_{1}+D_{2,n}}(m) = R_{H_{1}}(m) - 1$. 
For $u > 0$ even,  
to choose the pairs $(\varepsilon_{0}, \{z_{1}, \dots{} ,z_{w}\})$ we 
have $2\cdot 2^{u-1} = 2^{u}$ possibilities, thus we have $r_{D_{1}+D_{2,n}}(m) = R_{H_{1}}(m) - 2^{u}$. As $R_{H_{0}}(m) = R_{H_{1}}(m)$, it follows that 
$r_{C_{1}+C_{2,n}}(m) = r_{D_{1}+D_{2,n}}(m)$.

\noindent \textbf{Case 3.} $\delta_{0} = 2$. Then, $2 + u + 2v + t + 1$ is even, thus $u$ is odd. 
If $m = h + h^{'}$ with $h < h^{'}$ and $h, h^{'} \in H_{0}$, then $h \in C_{1}$ and $h^{'} \in H_{0} \setminus C_{1}$. It is clear that 
$h^{'} \notin C_{2,n}$ if and only if $h^{'}$ can be written in the form
\[
h^{'} = d_{1} + \sum_{j=1}^{w}d_{2^{z_{j}-1}+1} + \sum_{j=1}^{v}d_{2^{y_{j}-1}+1}  
+ d_{2^{i_{1}-1}+1} + d_{2^{i_{2}-1}+1} + \dots{} + d_{2^{i_{t}-1}+1} + d_{2^{N-1}+1},
\]
where $\{z_{1}, \dots{} ,z_{w}\} \subset \{x_{1}, \dots{} ,x_{u}\} \ne \emptyset$
and $1 + w + v + t + 1$ is even. There are $2^{u-1}$ ways to choose the set $\{z_{1}, \dots{} ,z_{w}\}$, thus we have $r_{C_{1}+C_{2,n}}(m) = R_{H_{0}}(m) - 2^{u-1}$.
Moreover, if $m = h + h^{'}$ with $h < h^{'}$ and $h, h^{'} \in H_{1}$, then  $h \in D_{1}$, $h^{'} \in H_{1} \setminus D_{1}$. It is clear that 
$h^{'} \notin D_{2,n}$ if and only if $h^{'}$ can be written in the form
\[
h^{'} = d_{1} + \sum_{j=1}^{w}d_{2^{z_{j}-1}+1} + \sum_{j=1}^{v}d_{2^{y_{j}-1}+1}  
+ d_{2^{i_{1}-1}+1} + d_{2^{i_{2}-1}+1} + \dots{} + d_{2^{i_{t}-1}+1} + d_{2^{N-1}+1},
\]
where $\{z_{1}, \dots{} ,z_{w}\} \subset \{x_{1}, \dots{} ,x_{u}\} \ne \emptyset$
and $1 + w + v + t + 1$ is odd. There are $2^{u-1}$ ways to choose the set $\{z_{1}, \dots{} ,z_{w}\}$, thus we have 
$R_{D_{1}+D_{2,n}}(m) = R_{H_{1}}(m) - 2^{u-1}$. As $R_{H_{0}}(m) = R_{H_{1}}(m)$, it follows that $r_{C_{1}+C_{2,n}}(m) = r_{D_{1}+D_{2,n}}(m)$.
If 
\[
K =  2d_{2^{i_{1}-1}+1} + d_{2^{i_{2}-1}+1} + \dots{} + d_{2^{i_{t}-1}+1} + d_{2^{N-1}+1} =  h + h^{'}
\]
with $h < h^{'}$ and $h, h^{'} \in H_{0}$, then
it follows from (5) and (10) that $h \in C_{1}$, $h^{'} \in H_{0} \setminus C_{1}$ and $h^{'}$ can be written in the form
\[
h^{'} = d_{2^{i_{1}-1}+1} + \sum_{j=1}^{w}d_{2^{z_{j}-1}+1} + d_{2^{N-1}+1},
\]
where  $\{z_{1}, \dots{} ,z_{w}\} \subset \{i_{2}, \dots{} ,i_{t}\} \ne \emptyset$. Thus, we have 
\[
h^{'} \le d_{2^{i_{1}-1}+1} + d_{2^{i_{2}-1}+1} + \dots{} + d_{2^{i_{t}-1}+1} + d_{2^{N-1}+1}
\]
\[
< d_{1} + d_{2^{i_{1}-1}+1} + d_{2^{i_{2}-1}+1} + \dots{} + d_{2^{i_{t}-1}+1} + d_{2^{N-1}+1} = H_{0,n+1},
\]
which implies $h^{'} \in C_{2,n}$ and $R_{H_{0}}(K) = r_{C_{1}+C_{2,n}}(K)$. 
In the last step, we prove $r_{C_{1}+C_{2,n}}(K) > r_{D_{1}+D_{2,n}}(K)$. It is clear that 
\[
K = d_{2^{i_{1}-1}+1} + ( d_{2^{i_{1}-1}+1} + d_{2^{i_{2}-1}+1} + \dots{} + d_{2^{i_{t}-1}+1} + d_{2^{N-1}+1}) =  d_{2^{i_{1}-1}+1} + H_{1,n+1},
\]
where $d_{2^{i_{1}-1}+1}, H_{1,n+1} \in H_{1}$. Since $H_{1,n+1} \notin D_{2,n}$, we have $R_{H_{1}}(K) > r_{D_{1}+D_{2,n}}(K)$. 
It follows from $R_{H_{0}}(K) = R_{H_{1}}(K)$ that $r_{D_{1}+D_{2,n}}(K) < r_{C_{1}+C_{2,n}}(K)$.

We will prove (iii) and (iv) simultaneously. Let 
\[
L = d_{2^{i_{1}-1}+1} + d_{2^{i_{2}-1}+1} + \dots{} + d_{2^{i_{t}-1}+1} + 2d_{2^{N-1}+1}. 
\]
We have to prove that if $m < L$, then $R_{D_{2,n}}(m) = R_{C_{2,n}}(m)$ and $R_{D_{2,n}}(L) < R_{C_{2,n}}(L)$. If $m < 2d_{2^{N-1}+1}$, then by using (7) and (11), we get that
$R_{D_{2,n}}(m) = R_{C_{2,n}}(m) = 0$. Assume that 
\[
2d_{2^{N-1}+1} \le m < L = d_{2^{i_{1}-1}+1} + d_{2^{i_{2}-1}+1} + \dots{} + d_{2^{i_{t}-1}+1} + 2d_{2^{N-1}+1}. 
\]
If $m = h + h^{'}$ with $h < h^{'}$ and $h, h^{'} \in H_{0}$, then it follows from (5) and (10) that $h, h^{'} \in H_{0} \setminus C_{1}$.
This implies that $h \ge d_{2^{N-1}+1}$ and 
\[
h^{'} < d_{2^{i_{1}-1}+1} + d_{2^{i_{2}-1}+1} + \dots{} + d_{2^{i_{t}-1}+1} + d_{2^{N-1}+1} < H_{0,n+1}.
\] 
It follows that $h, h^{'} \in C_{2,n}$, which yields $R_{H_{0}}(m) = R_{C_{2,n}}(m)$.
If $m = h + h^{'}$ with $h < h^{'}$ and $h, h^{'} \in H_{1}$, then it follows from (5) and (6) that $h, h^{'} \in H_{1} \setminus D_{1}$.
Since $h \ge d_{2^{N-1}+1}$ and 
\[
h^{'} < d_{2^{i_{1}-1}+1} + d_{2^{i_{2}-1}+1} + \dots{} + d_{2^{i_{t}-1}+1} + d_{2^{N-1}+1} = H_{1,n+1},
\]
it follows that $h, h^{'} \in D_{2,n}$, which yields $R_{H_{1}}(m) = R_{D_{2,n}}(m)$. As $R_{H_{0}}(m) = R_{H_{1}}(m)$, it follows that 
$R_{C_{2,n}}(m) = R_{D_{2,n}}(m)$. If 
\[
L = d_{2^{i_{1}-1}+1} + d_{2^{i_{2}-1}+1} + \dots{} + d_{2^{i_{t}-1}+1} + 2d_{2^{N-1}+1} = h + h^{'} 
\]
with $h < h^{'}$ and $h, h^{'} \in H_{0}$,
 then it follows from (5) and (10) that $h, h^{'} \in H_{0} \setminus C_{1}$. It follows that $h > d_{2^{N-1}+1}$ and
\[
h^{'} < d_{2^{i_{1}-1}+1} + d_{2^{z_{1}-1}+1} + \dots{} + d_{2^{z_{w}-1}+1} + d_{2^{N-1}} < H_{0,n+1},
\]
thus we have $h, h^{'} \in C_{2,n}$, which implies that $R_{H_{0}}(L) = R_{C_{2,n}}(L)$. 
On the other hand,
\[
L = d_{2^{i_{1}-1}+1} + d_{2^{i_{2}-1}+1} + \dots{} + d_{2^{i_{t}-1}+1} + 2d_{2^{N-1}+1} 
\]
\[
= d_{2^{N-1}+1} + (d_{2^{i_{1}-1}+1} + d_{2^{i_{2}-1}+1} + \dots{} + d_{2^{i_{t}-1}+1} + d_{2^{N-1}+1}) 
\]
\[
= d_{2^{N-1}+1} + H_{1,n+1}.
\]
Note that $H_{1,n+1}, d_{2^{N-1}+1} \in H_{1}$ and $H_{1,n+1} \notin D_{2,n}$, which
gives $R_{H_{1}}(L) > R_{C_{2,n}}(L)$. It follows from $R_{H_{0}}(L) = R_{H_{1}}(L)$
that $R_{D_{2,n}}(L) > R_{C_{2,n}}(L)$. 
\end{proof}

Now we are ready to prove that $H^{(n)}_{0}=C_1\cup C_{2,n}$ and $H^{(n)}_{1}=D_1\cup D_{2,n}$ hold for every $1 \le n \le 2^{N-1}$.
We prove by induction on $n$ that $C_{1} \cup C_{2,n} = H_{0}^{(n)}$ and 
 $D_{1} \cup D_{2,n} = H_{1}^{(n)}$. We have already proved $C_{1} \cup C_{2,1} = H_{0}^{(1)}$ and
 $D_{1} \cup D_{2,1} = H_{1}^{(1)}$.

\noindent Assume that $C_{1} \cup C_{2,n} = H_{0}^{(n)}$ and 
 $D_{1} \cup D_{2,n} = H_{1}^{(n)}$ hold for some $1 \le n < 2^{N-1}$. We will prove that $C_{1} \cup C_{2,n+1} = H_{0}^{(n+1)}$ and $D_{1} \cup D_{2,n+1} = H_{1}^{(n+1)}$ hold, i.e., $c_{2^{N-1}+n+1} = H_{0,n+1}$ and $d_{2^{N-1}+n+1} = H_{1,n+1}$.
Let 
\[
H_{0,n+1} = \varepsilon_{0}d_{1} + \sum_{j=1}^{t}d_{2^{i_{j}-1}+1} + d_{2^{N-1}+1},
\]
where $\varepsilon_{0} \in \{0, 1\}$, $(1 \le i_{1} < \dots{} < i_{t} < N)$.

\noindent \textbf{Case 1.} $\varepsilon_{0} = 0$, $t = 1$. We know from Lemma 1 that 
\[
min\{t_{k}, c_{2^{N-1}+1} + c_{2^{N-1}+n+1}\} = min\{s_{l}, d_{2^{N-1}+1} + d_{2^{N-1}+n+1}\}
\] 
and from Lemma 2 that $t_{k} = s_{l} = +\infty$. These facts imply that $c_{2^{N-1}+1} + c_{2^{N-1}+n+1} = d_{2^{N-1}+1} + d_{2^{N-1}+n+1}$.
Furthermore, we know that $c_{2^{N-1}+1} = d_{2^{N-1}+1} + d_{1}$, thus we have $c_{2^{N-1}+n+1} + d_{1} = d_{2^{N-1}+n+1}$, and then $d_{2^{N-1}+n+1} > c_{2^{N-1}+n+1}$. It follows from Lemma 1 that $min\{p_{i}, c_{2^{N-1}+n+1}\} = min\{q_{j}, d_{2^{N-1}+n+1}\}$ and from Lemma 2 that $p_{i} > q_{j} = H_{0,n+1}$. Then, we have $c_{2^{N-1}+n+1} = q_{j} =  H_{0 ,n+1}$ and   
\[
d_{2^{N-1}+n+1} = c_{2^{N-1}+n+1} + d_{1} = H_{0 ,n+1} + d_{1} = H_{1 ,n+1}.
\]
\noindent \textbf{Case 2.} $\varepsilon_{0} = 0$, $t > 1$. Applying Lemma 2, we get that $p_{i} > q_{j}$,
thus from Lemma 1 we have $min\{q_{j}, d_{1} + d_{2^{N-1}+n+1}\} = min\{p_{i}, c_{2^{N-1}+n+1}\} = c_{2^{N-1}+n+1}$. On the other hand, it follows from Lemma 2 that
$s_{l} < t_{k}$, thus by Lemma 1, we have 
\[
min\{s_{l}, d_{2^{N-1}+1} + d_{2^{N-1}+n+1}\} = min\{t_{k}, c_{2^{N-1}+1} + c_{2^{N-1}+n+1}\} 
\]
\[
= c_{2^{N-1}+1} + c_{2^{N-1}+n+1}. 
\]
Assume that $c_{2^{N-1}+n+1} = d_{1} + d_{2^{N-1}+n+1}$. Then, we have 
\begin{eqnarray*}
&&c_{2^{N-1}+1} + c_{2^{N-1}+n+1} = d_{1} + d_{2^{N-1}+1} + d_{1} + d_{2^{N-1}+n+1} \\
&=& 2d_{1} + d_{2^{N-1}+1} + d_{2^{N-1}+n+1} \\ 
&>& d_{2^{N-1}+1} + d_{2^{N-1}+n+1} \ge min\{s_{l}, d_{2^{N-1}+1} + d_{2^{N-1}+n+1}\}, 
\end{eqnarray*}
which is a contradiction. It follows from Lemma 2 that $c_{2^{N-1}+n+1} = q_{j} = H_{0,n+1}$ and 
\[
c_{2^{N-1}+1} + c_{2^{N-1}+n+1} = d_{1} + d_{2^{N-1}+1} + d_{2^{i_{1}-1}+1} + d_{2^{i_{2}-1}+1} + \dots{} + d_{2^{i_{t}-1}+1} + d_{2^{N-1}+1} 
\]
\[
= H_{1,n+1} + d_{2^{N-1}+1} = min\{s_{l}, d_{2^{N-1}+1} + d_{2^{N-1}+n+1}\} 
\]
\[
= min\{2d_{2^{i_{1}-1}+1} + d_{2^{i_{2}-1}+1} + \dots{} + d_{2^{i_{t}-1}+1} + 2d_{2^{N-1}+1}, d_{2^{N-1}+1} + d_{2^{N-1}+n+1}\}. 
\]
Since 
\[
d_{1} + d_{2^{N-1}+1} + d_{2^{i_{1}-1}+1} + \dots{} + d_{2^{i_{t}-1}+1} + d_{2^{N-1}+1} 
\]
\[
< 2d_{2^{i_{1}-1}+1} + \dots{} + d_{2^{i_{t}-1}+1} + 2d_{2^{N-1}+1}, 
\]
it follows that $d_{2^{N-1}+1} + d_{2^{N-1}+n+1} = H_{1,n+1} + d_{2^{N-1}+1}$, thus we have $d_{2^{N-1}+n+1} = H_{1,n+1}$.

\noindent \textbf{Case 3.} $\varepsilon_{0} = 1$. Applying Lemma 3, we get that $q_{j} > p_{i}$,
thus from Lemma 1, we have $min\{p_{i}, c_{2^{N-1}+n+1}\} = min\{q_{j}, d_{1} + d_{2^{N-1}+n+1}\} = d_{1} + d_{2^{N-1}+n+1}$. On the other hand, it follows from Lemma 3 that
$s_{l} > t_{k}$, thus by Lemma 1, we have 
\[
min\{t_{k}, c_{2^{N-1}+1} + c_{2^{N-1}+n+1}\} = min\{s_{l}, d_{2^{N-1}+1} + d_{2^{N-1}+n+1}\} 
\]
\[
= d_{2^{N-1}+1} + d_{2^{N-1}+n+1}. 
\]
Assume that $c_{2^{N-1}+1} + c_{2^{N-1}+n+1} = d_{2^{N-1}+1} + d_{2^{N-1}+n+1}$. 
Then, we have 
$d_{2^{N-1}+1} + d_{2^{N-1}+n+1} = d_{1} + d_{2^{N-1}+1} + c_{2^{N-1}+n+1}$, thus we have
$d_{2^{N-1}+n+1} = d_{1} + c_{2^{N-1}+n+1}$. It follows that
 $d_{1} + d_{2^{N-1}+1} = 2d_{1} + c_{2^{N-1}+n+1} = min\{p_{i}, c_{2^{N-1}+n+1}\}$, which is a contradiction because $d_{1} > 0$. Then, we have 
\[
d_{2^{N-1}+1} + d_{2^{N-1}+n+1} = t_{k} = d_{2^{i_{1}-1}+1} + d_{2^{i_{2}-1}+1} + \dots{} + d_{2^{i_{t}-1}+1} + 2d_{2^{N-1}+1}. 
\]
It follows that 
\[
d_{2^{N-1}+n+1} = d_{2^{i_{1}-1}+1} + d_{2^{i_{2}-1}+1} + \dots{} + d_{2^{i_{t}-1}+1} + d_{2^{N-1}+1} = H_{1,n+1}.
\]
Applying Lemma 1 and Lemma 3, we get that 
\begin{eqnarray*}
&&d_{1} + d_{2^{N-1}+n+1} = H_{0,n+1} \\ 
&=& d_{1} + d_{2^{i_{1}-1}+1} + d_{2^{i_{2}-1}+1} + \dots{} + d_{2^{i_{t}-1}+1} + d_{2^{N-1}+1}\\ 
&=& min\{p_{i}, c_{2^{N-1}+n+1}\}\\ 
&=& min\{2d_{2^{i_{1}-1}+1} + d_{2^{i_{2}-1}+1} + \dots{} + d_{2^{i_{t}-1}+1} + d_{2^{N-1}+1}, c_{2^{N-1}+n+1}\}\\ 
&=& c_{2^{N-1}+n+1}, 
\end{eqnarray*}
thus we have $c_{2^{N-1}+n+1} = H_{0,n+1}$. The proof of Theorem 4 has been completed.

\section{Proof of Theorem 6.}

First we prove that for 
\[
H=H(1,2,4,8,\dots, 2^{2l-1},2^{2l}-1,2^{2l+1}-1,2(2^{2l+1}-1),4(2^{2l+1}-1),8(2^{2l+1}-1),\dots ), 
\]
$C=H_0$ and $D=H_1$, we have $C\cup D=\mathbb{N}$, $C\cap D=2^{2l}-1+(2^{2l+1}-1)\mathbb{N}$ and $R_C(m)=R_D(m)$. It is easy to see that for $H^{'}=H(h_1,h_2,\dots ,h_{2l+1})=H(1,2,4,8,\dots ,2^{2l-1},2^{2l}-1)$, $C^{'}=H_0^{'}$ and $D^{'}=H_1^{'}$, we have $C^{'}\cup D^{'}=[0,2^{2l+1}-2]$ and $C^{'}\cap D^{'}=\{ 2^{2l}-1 \}$ because 
\[
2^{2l}-1=h_{2l+1}=h_1+h_2+\dots +h_{2l} = 1+2+4+\dots +2^{2l-1}. 
\]
Furthermore, for $H^{''}=H(2^{2l+1}-1,2(2^{2l+1}-1),4(2^{2l+1}-1),8(2^{2l+1}-1),\dots )$,

\noindent $C^{''}=H_0^{''}$ and $D^{''}=H_1^{''}$, we have $C^{''}\cup D^{''}=(2^{2l+1}-1)\mathbb{N}$ and $C^{''}\cap D^{''}=\emptyset$, which implies $C\cup D=\mathbb{N}$ and $C\cap D=2^{2l}-1+(2^{2l+1}-1)\mathbb{N}$. Moreover, by Theorem 3, $R_C(m)=R_D(m)$ for every positive integer $m$.
On the other hand, let us suppose that for some sets $C$ and $D$, we have $C\cup D=\mathbb{N}$ and $C\cap D=r+m\mathbb{N}$. By Conjecture 2, we may assume that for some Hilbert cube $H(h_1,h_2,\dots )$, we have $C=H_0$ and $D=H_1$. We have to prove the existence of integer $l$ such that $h_i=2^{i-1}$ for $1\le i\le 2l$, $h_{2l+1}=2^{2l}-1$ and $h_{2l+2+j}=2^{j}(2^{2l+1}-1)$ for every nonnegative integer $j$. We may suppose that $h_1=1$ and $h_2=2$.  Consider the Hilbert cube $H(1, 2, 4, \dots{}, 2^{u}, h_{u+2}, \dots{})$, where $h_{u+2}\ne 2^{u+1}$. Let us write $v=h_{u+2}$. We will prove that $v = 2^{u+1} - 1$.
Assume that $v > 2^{u+1}$. Then, it is clear that $2^{u+1} \notin H$ because
$1 + 2 + \dots{} + 2^{u} = 2^{u+1} - 1 < 2^{u+1}$. Thus, we have $v < 2^{u+1}$
i.e., $v \le 2^{u+1} - 1$.  Assume that
$v \le 2^{u+1} - 2$. Considering $v$ as a one term sum, it follows that $v \in
D$. Moreover, if $v = \sum_{i=0}^{u}\lambda_{i}2^{i}$, $\lambda _i\in \{ 0,1\}$, then $\sum_{i=0}^{u}\lambda_{i}$ must be even; otherwise, $v$ would have two different representations from $D$. It follows that $v \in C$ and $v + 1 =h_1+h_{u+2}\in C$. Furthermore, if we have a representation $v + 1 = \sum_{i=0}^{u}\delta_{i}2^{i}$, $\delta _i\in \{ 0,1\}$, then  $\sum_{i=0}^{u}\delta_{i}$ must be odd; otherwise, $v$ would have two different representations from $C$. This implies that $v + 1 \in D$, thus we have $v$,
$v + 1 \in C \cap D$. It follows that $C \cap D = \{v, v + 1, \dots{}\}$ is an arithmetic progression with common difference $1$. This implies that the generating
functions of the sets $C$ and $D$ are of the form
\[
C(z) = p(z) + \frac{z^{v}}{1 - z},
\]
where $p(z)$ is a polynomial, and
\[
D(z) = q(z) + \frac{z^{v}}{1 - z},
\]
where $q(z)$ is a polynomial, and
\[
p(z) + q(z) = 1 + z + z^{2} + \dots{} + z^{v-1} = \frac{1 - z^{v}}{1 - z}.
\]
 Since $R_{C}(n) = R_{D}(n)$, we have
$C^{2}(z) - D^{2}(z) = C(z^{2}) - D(z^{2})$. It follows that
\[
\Big(p(z) + \frac{z^{v}}{1 - z}\Big)^{2} - \Big(q(z) + \frac{z^{v}}{1 - z}\Big)^{2} = p(z^{2}) + \frac{z^{2v}}{1 - z^{2}} - q(z^{2})
- \frac{z^{2v}}{1 - z^{2}},
\]
which implies 
\[
p^{2}(z) - q^{2}(z) + \frac{2z^{v}}{1 - z}(p(z) - q(z)) = p(z^{2}) - q(z^{2}).
\]
Thus we have
\[
(p(z) - q(z)) \cdot \frac{1 + z^{v}}{1 - z} = p(z^{2}) - q(z^{2}).
\]
We get
\[
(p(z) - q(z)) \cdot (1 + z^{v}) = (p(z^{2}) - q(z^{2})) \cdot (1 - z).
\]
The leading coefficient in one side is $-1$ and the other side is $1$, which is a contradiction.
Thus we get that $v = 2^{u+1} - 1$. It follows that the Hilbert cube is of the
form $H(1, 2, 4, 8, \dots{}, 2^{u}, 2^{u+1} - 1, \dots{})$.
As $h_{u+2}=2^{u+1} - 1 = 1 + 2 + \dots{} + 2^{u}=h_1+h_2+\dots +h_{u+1}$, and $2^{u+1} - 1$ is regarded as a one-term
sum contained in $D$,
we have $u + 1$ must be even, i.e., $u + 1 = 2l$. It follows that there exists an integer $l$ such that $h_i=2^{i-1}$ for $1\le i\le 2l$ and $h_{2l+1}=2^{2l}-1$. It follows that $2^{2l} - 1 \in C \cap D$ and 
$r = 2^{2l} - 1$.  

We apply induction on $j$ to show that $h_{2l+2+j} = 2^{j}(2^{2l+1} - 1)$ for every nonnegative integer $j$. For $j=0$, take the
Hilbert cube of the form
$H(1, 2, 4, 8, \dots{}, 2^{2l-1}, 2^{2l} - 1, h_{2l+2}, \dots{})$. 
Write $w=h_{2l+2}$. We prove that $w=2^{2l+1}-1$. Suppose that
$w > 2^{2l+1} - 1$.
Since $1 + 2 + \dots{} + 2^{2l-1} + 2^{2l} - 1 < 2^{2l+1} - 1$,
it follows that $2^{2l+1} - 1 \notin H = C \cup D$, which is impossible. Therefore,
$w \le 2^{2l+1} - 1$. Suppose that $w \le 2^{2l+1} - 3$. We will show that
$w \in C \cap D$. Obviously, $w$ is a one-term sum contained in $D$.
Since $w$ has a representation from $H(h_1,\dots ,h_{2l+1})$, $w$ must be an element of $C$; otherwise, $w$ would have two different
representations from $D$, which is absurd. In the next step, we prove
$w + 1 \in C \cap D$. Obviously, $w + 1=h_1+h_{2l+2}$ as a two-term sum contained in $C$.
Since $w +1$ can be written in terms of the Hilbert cube $H(h_1,\dots ,h_{2l+1})$ and $w + 1 \le 2^{2l+1} - 2$, we have $w + 1 \in D$. It follows that $w$, $w + 1 \in  C \cap D$, which is
a contradiction. It follows that the only possible values of $w$ are $w = 2^{2l+1} - 2$, and $w = 2^{2l+1} - 1$. Suppose that $w = 2^{2l+1} - 2$. Then, it is clear that $w \in D$. On the other hand,
$2^{2l} - 2 = 1 + 2 + \dots{} + 2^{2l-1} + 2^{2l} - 1=h_1+h_2+\dots +h_{2l+1}$, where in the right hand side, there are $2l + 1$ terms, which is impossible. Thus we have $w = 2^{2l+1} - 1$. In this case, $2^{2l} - 1$, $(2^{2l}-1)+(2^{2l+1}-1)\in C\cap D$, $(C\cap D)\cap \{1,2,\dots ,2^{2l+1}-1\}=\{2^{2l}-1\}$. It follows that $m \mid 2^{2l+1}-1$. If $m \le \frac{2^{2l+1}-1}{2}$, then $(C \cap D) \cap \{1, 2, \dots{}, 2^{2l+1}-1\} \ne \{2^{2l}-1\}$, a contradiction. Then, we have $r = 2^{2l}-1$ and $m =  2^{2l+1}-1$.

In the induction step, we assume that for some $k$, we know that $h_{2l+2+j}=2^{j}(2^{2l+1}-1)$ holds for $j=0,1,\dots ,k$, and we prove $h_{2l+2+k+1}=2^{k+1}(2^{2l+1}-1)$. Let 
$H^{(k)}=H(1,2,4,8,\dots ,2^{2l-1},2^{2l}-1,2^{2l+1}-1,2(2^{2l+1}-1),4(2^{2l+1}-1),8(2^{2l+1}-1),\dots ,2^{k}(2^{2l+1}-1))$, 
$C^{(k)}=H_0^{(k)}$ and $D^{(k)}=H_0^{(k)}$. Then 
\[
C^{(k)}\cap D^{(k)}=\{2^{2l}-1+i(2^{2l+1}-1): i=0,1,\dots ,2^{k}-1\}. 
\]
If $C=H_0(1,2,4,8,\dots ,2^{2l-1},2^{2l}-1,2^{2l+1}-1,2(2^{2l+1}-1),4(2^{2l+1}-1),8(2^{2l+1}-1),\dots ,2^{k}(2^{2l+1}-1),h_{2l+2+k+1},\dots )$, 
then $C\cap D=\{ e_1,e_2,\dots \}$, where $e_i=2^{2l}-1+(i-1)(2^{2l+1}-1)$ for $i \ge 1$, and $e_{2^{k+1}+1}=2^{2l}-1+2^{k+1}(2^{2l+1}-1)$. Furthermore, $e_{2^{k+1}+1}=2^{2l}-1+h_{2l+1+k+1}$, and then $h_{2l+1+k+1}=2^{k+1}(2^{2l+1}-1)$, which completes the proof.



\begin{thebibliography}{99}
\bibitem{CH} Y. G. Chen, On the values of representation
  functions, {\it Sci. China Math.} {\bf 54} (2011), 1317-1331.
\bibitem{CL} Y. G. Chen and V. F. Lev, Integer sets with
  identical representation functions, {\it Integers} {\bf 16} (2016), A36.
\bibitem{CT} Y. G. Chen and M. Tang, Partitions of natural
  numbers with the same representation functions, {\it J. Number Theory}
{\bf 129} (2009), 2689-2695.
\bibitem{CW} Y. G. Chen and B. Wang, On additive properties of
  two special sequences, {\it Acta Arith.}, {\bf 113} (2003), 299-303.
\bibitem{DG} G. Dombi, Additive properties of certain sets,
  {\it Acta Arith.} {\bf 103} (2002), 137-146.
\bibitem{KR} S.Z. Kiss and E. Rozgonyi and Cs. S\'andor,
Sets with almost coinciding representation functions, {\it Bull. Aust. Math. Soc.}, {\bf 89} (2014), 97-111.
\bibitem{FL} V. F. Lev, Reconstructing integer sets from
  their representation functions, {\it Electron. J. Combin.}, {\bf 11} (2004),
R78.
\bibitem{NT} M. B. Nathanson, Representation
  functions of sequences in additive number theory, {\it Proc.
Amer. Math. Soc.}, {\bf 72} (1978), 16-20.
\bibitem{ZQ} Z. Qu, On the nonvanishing of representation
  functions of some special sequences, {\it Discrete Math.}, {\bf 338} (2015),
  571-575.
\bibitem{RS} E. Rozgonyi and Cs. S\'andor, An extension of Nathanson's Theorem on representation functions,
{\it Combinatorica}, {\bf 164} (2016), 1-17.
\bibitem{CS} Cs. S\'andor, Partitions of natural numbers and
  their representation functions, {\it Integers}, {\bf 4} (2004), A18.
\bibitem{SS} J. L. Selfridge and E. G. Straus, On the determniation of numbers by their sums of a fixed order, {\it Pacific J. Math.}, {\bf 8} (1958), 847-856.
\bibitem{MT} M. Tang, Partitions of the set of natural
numbers and their representation functions, {\it Discrete Math.}, {\bf 308}
(2008), 2614-2616.
\bibitem{MY} M. Tang and W. Yu, A note on partitions of natural
  numbers and their representation functions, {\it Integers}, {\bf 12} (2012),
A53.


\end{thebibliography}
\end{document}